\documentclass[11pt]{article}

\usepackage{amssymb}
\usepackage{amsmath}
\usepackage{latexsym}
\usepackage{mathrsfs}

\newtheorem{theorem}{Theorem}

\newtheorem{lemma}[theorem]{Lemma}

\newtheorem{definition}[theorem]{Definition}
\newtheorem{remark}[theorem]{Remark}

\newcommand{\po}{{\partial\Omega}}

\newcommand{\dom}{{\rm Dom}}

\newcommand{\kernel}{{\rm Ker}}
\newcommand{\ran}{{\rm Ran}}
\newcommand{\dist}{{\rm dist}}
\newcommand{\N}{{\bf N}}

\newcommand{\R}{{\bf R}}
\newcommand{\C}{{\bf C}}
\newcommand{\hil}{{\cal H}}
\newcommand{\cC}{{\cal C}}

\newcommand{\cD}{{\cal D}}

\newcommand{\cV}{{\cal V}}
\newcommand{\cW}{{\cal W}}

\newcommand{\parder}[2]{\frac{\partial{#1}}{\partial{#2}}}

\newcommand{\inprod}[2]{{\langle{#1},{#2}\rangle}}

\numberwithin{equation}{section}

\title{Stability estimates for resolvents, eigenvalues and eigenfunctions of elliptic operators on variable domains}
\author{G. Barbatis, V.~I. Burenkov, P.~D. Lamberti}
\date{}
\begin{document}
\maketitle

\begin{center}
{\em Dedicated to  Vladimir Maz'ya}
\end{center}

\begin{abstract}
We consider general second order uniformly elliptic operators
subject to homogeneous boundary conditions on open sets $\phi (\Omega )$ parametrized  by  Lipschitz homeomorphisms $\phi $ defined on a fixed reference domain $\Omega$.
Given two open sets $\phi (\Omega )$, $\tilde \phi (\Omega )$ we estimate  the variation of resolvents, eigenvalues and eigenfunctions
  via the Sobolev norm $\|\tilde  \phi -\phi \|_{W^{1,p}(\Omega )}$ for finite values of $p$, under natural summability conditions on eigenfunctions
  and their gradients. We prove that such conditions are satisfied for a wide  class of operators and open sets, including open sets with Lipschitz continuous boundaries.  We apply these estimates  to control the variation of the eigenvalues and eigenfunctions via the measure of the symmetric difference of the open sets. We also discuss  an application to the stability of solutions to the Poisson problem.
\end{abstract}

\vspace{11pt}


\vspace{6pt}
%



\

\section{Introduction}

This paper is devoted to the proof of stability estimates for the non-negative self-adjoint operator

\begin{equation}\label{intro1}
Lu=-\sum_{i,j=1}^N\frac{\partial}{\partial x_i }\big( A_{ij}(x)\frac{\partial u}{\partial x_j} \big),\ \ \ x\in\Omega ,
\end{equation}
subject to homogeneous boundary conditions, upon variation of the open set $\Omega $ in $\R^N$. Here $A_{ij}$ are fixed bounded measurable real-valued functions defined in $\R^N$ satisfying $A_{ij}=A_{ji}$ and a uniform ellipticity condition.

The focus is on explicit quantitative estimates for the variation of the resolvents,  eigenvalues and  eigenfunctions of $L$ on a class of open sets diffeomorphic to $\Omega$.

In the first part of the paper we  consider two diffeomorphisms $ \phi, \tilde \phi  $ of $\Omega $ onto $\phi (\Omega )$, $\tilde \phi (\Omega )$ respectively, and we  compare the resolvents, eigenvalues and eigenfunctions of $L$ on the open set  $\tilde\phi (\Omega )$ with those of $L$ on $\phi (\Omega )$.
The main aim is to provide stability estimates via $\|\tilde\phi- \phi\|_{W^{1,p}(\Omega)}$ for finite values of $p$.
These estimates are applied in the last part of the paper where we  take $\phi =Id$ and, given a deformation $\tilde \Omega$ of $\Omega$, we  construct a special diffeomorphism $\tilde \phi$ representing $\tilde \Omega$ in the form $\tilde \Omega =\tilde\phi (\Omega  )$  and obtain stability estimates in terms of the Lebesgue measure $|  \Omega\vartriangle \tilde \Omega|$ of the symmetric difference of $\Omega$ and $\tilde\Omega$.

 Our method allows us to treat the general case of the mixed homogeneous Dirichlet-Neumann boundary conditions
\begin{equation}
\label{intro2}
u=0 \ {\rm on}\ \Gamma ,\ \ \ {\rm and}\ \ \
\sum_{i,j=1}^N A_{ij}\frac{\partial u}{\partial x_j}\nu_{i}=0 \ {\rm on}\ \partial\Omega  \setminus \Gamma ,
\end{equation}
where $\Gamma\subset\partial\Omega$. To our knowledge, our results are new also for  Dirichlet, and for  Neumann boundary conditions.

There is a vast literature concerning domain perturbation problems, see for instance the extensive monograph by Henry~\cite{he}.
The problem of finding explicit quantitative estimates  for the variation of the eigenvalues
of elliptic operators has been considered in Burenkov and Davies~\cite{buda}, Burenkov and Lamberti~\cite{bulahigh, bula, bulaneu}, Burenkov and Lanza de Cristoroforis~\cite{bulanz}, Davies~\cite{da2000, da1993}, Lamberti and Lanza de Cristoforis~\cite{lala03,lala03neu}, and Pang~\cite{pang}; see  Burenkov, Lamberti and Lanza de Cristoforis~\cite{bulalanz} for a survey on the results of these papers. However, less attention has been devoted to the problem of finding explicit estimates for the variation of the eigenfunctions. With regard to this, we mention the estimate in \cite{pang} concerning the first eigenfunction of the Dirichlet Laplacian and the estimates in \cite{lala03, lala03neu} concerning the variation of the eigenprojectors of the Dirichlet and Neumann Laplacian. In particular, in \cite{lala03, lala03neu}  the variation of the eigenvalues and  eigenprojectors of the Laplace operator was estimated via $\| \nabla\tilde  \phi -\nabla\phi \|_{L^{\infty }(\Omega )}$
under minimal assumptions on the regularity of $\Omega$, $\phi$ and  $\tilde \phi $.

In all cited papers and in this paper perturbations of domains may be considred as in some sense regular pertubations. There is also vast literature concerning a wide range of perturbation problems of different type which may be characterised as singular perturbations (which are out of scope of this paper). Typically, formulations of such problems involve a small parameter $\varepsilon$ and the problem degenerates in that sense or other as $\varepsilon\to 0$. Say, the domain may contain small holes, or boundaries which may include blunted angles, cones and edges, narrow slipts, thin bridges etc, or the limit region may consist of subsets of different dimension, or it could be a homogenization problem. V.~G.~Maz'ya and his co-authors V.~A.~Kozlov, A.~B.~Movchan, S.~A.~Nazarov, B.~A.~Plamenevskii and others developed the powerful asymptotic theory which allowed to find asymptotic expansions of solutions for all aforementioned problems and can be applied in many other cases. See  {\it e.g.}, monographs \cite{m5, m1, m2}. 

In this paper, we consider the same class of transformations $\phi , \tilde \phi$ as in \cite{lala03, lala03neu} ($\phi , \tilde \phi $ are bi-Lipschitz homeomorphisms) and by making stronger regularity assumptions on $\phi (\Omega )$, $\tilde \phi (\Omega )$ we estimate  the variation of the resolvents,  eigenvalues, eigenprojectors and     eigenfunctions of $L$ via   the measure of vicinity
\begin{equation}
\label{deltaphi}
\delta_p (\phi ,\tilde\phi ):=\| \nabla \tilde \phi -\nabla\phi \|_{L^{p }(\Omega )}+\| A\circ\tilde \phi -A\circ  \phi\|_{L^{p }(\Omega )}
\end{equation}
for any $p\in ]p_0, \infty ]$, where $A=(A_{ij})_{i,j=1, \dots , N}$ is the matrix of the coefficients.

 Here $p_0\geq 2$ is a constant  depending on the regularity assumptions. The best $p_0$ that we obtain is $p_0=N$  which  corresponds
to the highest degree  of regularity (see Remark~\ref{bestrange}), whilst the case $p_0=\infty$ corresponds to the lowest degree of regularity in which case only the exponent  $p=\infty$ can be considered. The regularity assumptions are expressed in terms of summability properties of the eigenfunctions and their gradients, see Definition~\ref{growth}. Observe that if the coefficients $A_{ij}$ of the operator $L$  are Lipschitz continuous, then $\delta _{p}(\phi ,\tilde \phi)$ does not exceed a constant independent of $\phi,\tilde\phi$ multiplied by the Sobolev norm $\|\phi -\tilde \phi  \|_{W^{1,p}(\Omega )}$.
Moreover if the coefficients $A_{ij}$ are constant then the second summand in the right-hand side of (\ref{deltaphi}) vanishes.

More precisely, we prove stability estimates for the resolvents in the Schatten classes (Theorem~\ref{mainthm}), stability
estimates for eigenvalues (Theorem~\ref{thm:series}), eigenprojectors (Theorem~\ref{thm:projectors}), and eigenfunctions (Theorem~\ref{holder}). In the  Appendix we also consider an application to the Poisson problem (we refer to Savar\'{e} and Schimperna~\cite{sava} for stability estimates for the solutions to the Poisson problem in the case of Dirichlet boundary conditions obtained by a different approach).
To prove the  resolvent stability estimates  in the Schatten classes we follow the method developed in Barbatis \cite{bar1, bar3}.

In Section~\ref{applsec} we apply our general  results and, for a given deformation $\tilde\Omega$ of $\Omega$, we prove stability estimates in terms of $|\Omega\vartriangle\tilde \Omega |$.  This is done  in two cases: the case in which $\tilde \Omega $ is obtained by a localized deformation of the boundary of $\Omega $ and the case in which  $\tilde \Omega $ is a deformation of $\Omega $ along its normals. We also require that the deformation $\tilde\Gamma$ of $\Gamma  $ is induced by the deformation of $\Omega$ (see conditions (\ref{bdmixed}) and (\ref{tildegamma})).  In these cases, similarly to \cite{bula}, we can construct special bi-Lipschitz transformations $\tilde \phi:\Omega \to \tilde\Omega$ such that $\tilde \phi (\Gamma )=\tilde \Gamma $ and
\begin{equation}\label{trick}
\|\nabla\tilde \phi-I \|_{L^p(\Omega )}\le c|\Omega \vartriangle \tilde \Omega |^{1/p}.
\end{equation}
Observe  that using finite values of $p$ is essential, since in the case $p=\infty$ the exponent in the right-hand side of (\ref{trick}) vanishes.

Let us  describe these results  in the regular case in which $\Omega$, $\tilde\Omega$ are of class $C^{1,1}$ and $\Gamma $, $\tilde\Gamma$  are connected components of the corresponding boundaries.
In Theorems~\ref{thm:C11}, \ref{thm:C11bis}  we prove that for any $r>N$ there exists a constant $c_1 >0$ such that
\begin{equation}\label{introuno}
\left(\sum_{n=1}^{\infty}\left| \frac{1}{
{\lambda}_n +1} - \frac{1}{{\tilde \lambda}_n +1}\right|^r\right)^{1/r}
\leq c_1 |\Omega \vartriangle \tilde \Omega  |^{\frac{1}{r}} ,
\end{equation}
if $|\Omega \vartriangle\tilde\Omega  |<c_1^{-1}$. Here ${\lambda}_n$, ${\tilde \lambda}_n$ are the eigenvalues of the  operators (\ref{intro1}) corresponding to the domains $\Omega$, $\tilde\Omega$ and the associated portions of the boundaries $\Gamma $, $\tilde\Gamma$ respectively.
Moreover, for a fixed $\Omega$ and for any $r >N$ there exists $c_2>0$ such that  if $\lambda_n=\dots $ $=\lambda_{n+m-1}$ is an eigenvalue of multiplicity $m$  then for any choice
of orthonormal eigenfunctions
$\tilde \psi_n, \dots , $ $ \tilde \psi_{n+m-1}$ corresponding to $\tilde \lambda_n,\dots ,\tilde \lambda_{n+m-1}$, there exist  orthonormal eigenfunctions $ \psi_n,  \dots $, $    \psi_{n+m-1} $
corresponding to $\lambda_n, \dots ,$ $\lambda_{n+m-1}$ such that\footnote{Note that, for a fixed $\Omega$ and variable $\tilde\Omega $,  one first chooses eigenfunctions in $\tilde \Omega$ and then finds eigenfunctions in $\Omega$, while the opposite  is clearly not possible.}
\begin{equation}\label{introdue}
\| \psi_k-\tilde \psi_k \|_{L^2(\Omega \cup\tilde \Omega)}\le c_2 |\Omega \vartriangle\tilde \Omega |^{\frac{1}{r}},
\end{equation}
for all $k=n, \dots , n+m-1$, provided that $|\Omega \vartriangle\tilde\Omega  |<c_2^{-1}$. Here it is understood that the eigenfunctions are extended by zero outside their domains of definition.

In the general case of open sets $\Omega $, $\tilde\Omega$ with  Lipschitz continuous boundaries and  $\Gamma $, $\tilde\Gamma $ with Lipschitz continuous boundaries in $\partial \Omega $, $\partial \tilde\Omega $  our  statements still hold  for a possibly worse  range of exponents (see Theorems~\ref{thm:C11}, \ref{thm:C11bis}).

We emphasize that, in the spirit of \cite{lala03, lala03neu}, in this paper we never assume that the transformation $\phi $
belongs to a  family of transformations $\phi_{t}$ depending analytically on one scalar parameter  $t$, as often done in the literature (see {\it e.g.}, \cite{he} for references). In that case, one can use
 proper methods of bifurcation theory in order to prove existence of branches of eigenvalues and eigenfunctions depending analytically on $t$. In this paper $\tilde \phi$ is an arbitrary
perturbation of $\phi$  and this requires a totally different approach.

The paper is organized as follows: in Section~\ref{operators} we describe the general setting; in Section~\ref{perturb} we describe our perturbation problem; in Section~\ref{resolv} we prove
stability estimates for the resolvents and the eigenvalues;
in Section~\ref{eigenprsec} we prove stability estimates for the eigenprojectors and eigenfunctions;
in Section~\ref{reg} we give sufficient conditions providing the required regularity of the eigenfunctions;
in Section~\ref{applsec}  we prove stability estimates via  the Lebesgue
measure of the symmetric difference of sets; in the Appendix we briefly discuss the Poisson problem.

\

\section{General setting}
\label{operators}

Let $\Omega$ be a domain, {\it i.e.} an open connected set,  in $\R^N$ of finite measure. We
consider a family of
open sets $\phi\left(\Omega\right)$ in $\R^N$
parametrized by bi-Lipschitz homeomorphisms $\phi$ of $\Omega $ onto $\phi (\Omega )$.  Namely, following \cite{lala03}, we consider the  family of transformations \begin{eqnarray}
\label{phis}
& &\Phi (\Omega) := \left\{
\phi\in\left(L^{1,\infty}(\Omega)\right)^{N}:\,
{\mathrm{the\ continuous\ representative\ of}}\ \phi \right.\nonumber
\\
& &\left.
\qquad\qquad
\qquad\qquad\qquad\qquad\quad
{\mathrm{is\ injective}},\
{\mathrm{ess}}
\inf_{\Omega}|\det \nabla \phi |
>0
\right\},
\end{eqnarray}
where $L^{1,\infty}(\Omega)$ denotes the space
of the functions in $L^1_{loc}\left(\Omega\right)$
which have weak derivatives of first order in
$L^{\infty}\left(\Omega\right)$. Observe that if $\phi\in \Phi (\Omega )$ then $\phi$ is Lipschitz continuous with respect to the geodesic distance in $\Omega$.

Note that if $\phi \in \Phi (\Omega )$ then
$\phi (\Omega)$ is open, $\phi$ is a
homeomorphism of $\Omega$ onto $\phi(\Omega)$ and the inverse
vector-function $\phi^{(-1)}$ of $\phi$ belongs to
$\Phi\left(\phi (\Omega )\right)$.
Moreover, any transformation $\phi \in \Phi (\Omega )$ allows changing variables in integrals. Accordingly,
the operator
$C_{\phi}$ from $L^{2}(\phi (\Omega))$ to $L^{2}(\Omega)$
defined by
$$  C_{\phi}[v]:=
v\circ\phi,\ \  \ v\in
L^{2}(\phi (\Omega)),$$  is a linear homeomorphism
which restricts to a linear homeomorphism of the space
$W^{1,2}(\phi(\Omega))$ onto $W^{1,2}(\Omega)$, and of
$W^{1,2}_{0}(\phi(\Omega))$ onto $W^{1,2}_{0}(\Omega)$, where $W^{1,2}(\Omega) $ denotes the standard Sobolev space and $W^{1,2}_{0}(\Omega)$ denotes the closure of $C^{\infty }_c(\Omega )$ in $W^{1,2}(\Omega)$.
Furthermore, $\nabla (v\circ\phi )=\nabla v(\phi )\nabla \phi$ for all $v\in W^{1,2}\left(
\phi\left(\Omega\right)\right)$.
Observe that if $\phi\in \Phi (\Omega )$ then  the measure of $\phi (\Omega )$ is finite. See \cite{lala03} for details.

\

Let $A=(A_{ij})_{i,j=1, \dots , N}$ be a real symmetric matrix-valued measurable function defined on $\R^N$ such that for some $\theta >0$
\begin{equation}
\label{ellip}
\theta^{-1} |\xi |^2 \le \sum_{i,j=1}^NA_{ij}(x)\xi_i\xi_j \le \theta |\xi |^2 ,
\end{equation}
for all $x, \xi \in \R^N$. Note that (\ref{ellip}) implies that $A_{ij}\in L^{\infty }(\R^N )$ for all $i,j=1, \dots , N$.

Let $\phi\in \Phi(\Omega)$ and let ${\mathcal{W}}$ be a closed subspace of $W^{1,2}(\phi(\Omega))$
containing $W^{1,2}_0(\!\phi(\Omega))$. We consider the non-negative self-adjoint operator $L$ on $L^2(\phi(\Omega))$
canonically associated with the sesquilinear form $Q_L$ given by
\begin{equation}
\label{quadratic}
\dom(Q_L)={\mathcal{W}} \, , \quad Q_L(v_1,v_2)=\int_{\phi(\Omega)} \sum_{i,j=1}^NA_{ij}\parder{v_1}{y_i}\parder{\bar v_2}{y_j}
   dy \, , \ \ v_1,v_2\in {\mathcal{W}}.
\end{equation}
Recall that $v\in \dom(L)$ if and only if $v\in {\mathcal{W}}$ and there exists $f\in L^2(\phi (\Omega ))$ such that
\begin{equation}
\label{doml}
 Q_L(v,\psi )=\inprod{f}{\psi}_{L^2(\phi (\Omega ))},
\end{equation}
for all $\psi \in {\mathcal{W}}$, in which case $Lv=f$ (see {\it e.g.}, Davies~\cite{daheat}).
The choice of the space ${\mathcal{W}}$ determines the boundary conditions. For example if ${\mathcal{W}}=W^{1,2}_0(\phi (\Omega ))$
(respectively, ${\mathcal{W}}=W^{1,2}(\phi (\Omega ))$) then the operator $L$ satisfies homogeneous Dirichlet (respectively, homogeneous Neumann) boundary conditions.

We also consider the operator $H$ on $L^2(\Omega)$ obtained by pulling-back $L$ to $L^2(\Omega)$ as follows.
Let $v\in W^{1,2}(\phi (\Omega ))$ be given and let $u=v\circ \phi$. Observe that
\[
\int_{\phi(\Omega)}|v|^2dy =\int_{\Omega}|u|^2|\det\nabla\phi | \, dx\; .
\]
Moreover a simple computation shows that
\[
\int_{\phi(\Omega)}   \sum_{i,j=1}^NA_{ij}\parder{v}{y_i}\parder{\bar v}{y_j} dy =\int_{\Omega }\sum_{i,j=1}^Na_{ij}\parder{u}{x_i}\parder{\bar u}{x_j}
|\det\nabla\phi | \, dx\; ,
\]
where $a=(a_{ij})_{i,j=1, \dots , N}$ is the symmetric matrix-valued function defined on $\Omega$ by
\begin{eqnarray}
a_{ij}&=&\sum_{r,s=1}^N \Big(A_{rs}  \parder{  \phi_i^{(-1)}}{y_r}   \parder{\phi_j^{(-1)}}{y_s}\Big)\circ \phi \nonumber \\
&=& ((\nabla \phi)^{-1} A(\phi)  (\nabla\phi)^{-t})_{ij}\, .
\label{matrixa}
\end{eqnarray}
The operator $H$ is defined as the non-negative self-adjoint operator on the Hilbert space $L^2(\Omega , |\det\nabla\phi |\, dx)$
associated with the sesquilinear form $Q_H$ given by
\[
\dom(Q_H  ) =C_{\phi}[ {\mathcal{W}}  ],\
Q_{H}(u_1,u_2)=\int_{\Omega}\sum_{i,j=1}^Na_{ij}\parder{u_1}{x_i}\parder{\bar u_2}{x_j}|\det\nabla\phi |\, dx, \ u_1,u_2\in C_{\phi}[{\mathcal{W}}].
\]
Formally,
$$
Hu= -\frac{1}{|{\rm det}\nabla \phi | }\sum_{i,j=1}^N\frac{\partial }{\partial x_j}\left(a_{ij}\frac{\partial u}{\partial x_i} |{\rm det}\nabla \phi | \right) .
$$
Alternatively, the operator $H$ can be defined as
$$
H=C_{\phi }L C_{\phi ^{(-1)}}.
$$
In particular $H$ and $L$ are unitarily equivalent and the operator $H$ has compact resolvent if and only if $L$ has compact resolvent. (Observe that the embedding ${\mathcal{W}}  \subset L^2(\phi(\Omega))$ is compact if and only if the embedding $C_{\phi}[  {\mathcal{W}}]\subset L^2(\Omega)$ is  compact.)

We set $g(x):=|\det\nabla\phi(x)|$, $x\in\Omega$, and we denote by $\inprod{\cdot}{\cdot}_g$
the inner product in $L^2(\Omega,g\, dx)$ and also in $(L^2(\Omega,g\, dx))^N$.

Let $T:L^2(\Omega, g\, dx)\to (L^2(\Omega, g\, dx))^N$ be the operator defined by
$$\dom(T)=C_{\phi}[{\mathcal{W}}],\ \ \ \ \ Tu=a^{1/2}\nabla u,\ \ u\in C_{\phi}[{\mathcal{W}}].$$ It is then easily seen that
$$H=T^*T,$$
where the adjoint $T^*$ of $T$ is understood with respect to the
inner products of $L^2(\Omega, g\, dx)$ and $(L^2(\Omega, g\, dx))^N$.

\

\section{ Perturbation of $\phi $}
\label{perturb}

In this paper we study the variation of the operator $L$ defined by \eqref{quadratic} upon variation of $\phi$.
Our estimates depend on ${\mathrm{ess}}
\inf_{\Omega}|\det \nabla \phi |$ and $\|\nabla \phi  \|_{L^{\infty }(\Omega )}$. Thus in order to obtain uniform estimates it is convenient to consider the families of transformations
$$
\Phi_{\tau}(\Omega )=\left\{\phi\in \Phi (\Omega ):\ \tau ^{-1}\le {\mathrm{ess}}
\inf_{\Omega}|\det \nabla \phi |\ {\rm and }\  \|\nabla \phi  \|_{L^{\infty }(\Omega )} \le \tau    \right\},
$$
for all $\tau >0$, as  in \cite{lala03}. Here and below for a matrix-valued function $B(x)$, $x\in \Omega$,
we set  $\| B\|_{L^p(\Omega )}= \|\, |B |\, \|_{L^p(\Omega )}$ where $| B(x) |$ denotes the operator norm of the matrix $B(x)$.

Let $\phi,\tilde{\phi}\in \Phi_{\tau }(\Omega)$. Let ${\mathcal{W}}$ and $\tilde{{\mathcal{W}}}$ be closed subspaces of $W^{1,2}(\phi (\Omega ))$, $W^{1,2}(\tilde \phi (\Omega ))$ respectively, containing $W^{1,2}_0(\phi (\Omega ))$, $W^{1,2}_0(\tilde \phi (\Omega ))$ respectively.
We use tildes to distinguish  objects induced by $\tilde \phi$, $\tilde {\mathcal{W}}$  from  those induced by ${\phi}$, $ {\mathcal{W}}$.
We consider the operators $L$ and $\tilde L$  defined on $L^2(\phi (\Omega ))$, $L^2(\tilde \phi (\Omega ))$ respectively, as in Section~\ref{operators}.

In order to compare $L$ and $\tilde L$ we shall make a `compatibility' assumption on the respective boundary conditions; namely, we shall assume
that
\begin{equation}\label{compa} C_{\phi }[{\mathcal{W}}] = C_{\tilde \phi }[\tilde {\mathcal{W}}].
\end{equation}
This means that $\dom(Q_H)=\dom(Q_{\tilde H})$, a property which is important in what follows.
Clearly (\ref{compa}) holds  if either $L$ and $\tilde L$ both satisfy homogeneous Dirichlet boundary conditions or they both satisfy homogeneous Neumann
boundary conditions.

We shall always assume that the spaces ${\mathcal{W}}$, $\tilde {\mathcal{W}}$ are compactly embedded in $L^2(\phi (\Omega ))$, $L^2(\tilde \phi (\Omega ))$ respectively, or equivalently that the space ${\mathcal {V}}:=  C_{\phi }[{\mathcal{W}}] = C_{\tilde \phi }[\tilde {\mathcal{W}}]  $ is compactly embedded in $L^2(\Omega)$.

Moreover, we  require  that the non-zero eigenvalues $\lambda_n$ of the Laplace operator associated in $L^2(\Omega )$ with the quadratic form
$\int_{\Omega }|\nabla u|^2dx  $, $u\in {\mathcal {V}}$, defined on ${\mathcal{V}}$, satisfy the condition
\begin{equation}
c^*:=\sum_{\lambda_n\ne 0}\lambda_n^{-\alpha }<\infty\, ,
\label{cstar}
\end{equation}
for some fixed $ \alpha>0$. (This is in fact a very weak condition on the regularity of the set $\Omega $ and the associated boundary conditions.)

For brevity, we shall refer to assumption {\rm (A)} as the following set of conditions  which summarize the setting described above:

$${\rm (A)}:\ \ \left\{\begin{array}{l}
\phi, \tilde\phi \in \Phi_{\tau }(\Omega ),\vspace{1mm}\\
{\mathcal {V}}:= C_{\phi }[{\mathcal{W}}] = C_{\tilde \phi }[\tilde {\mathcal{W}}]\ {\rm is \ compactly \ embedded \ in }\ L^2(\Omega ),\vspace{1mm}\\
{\rm condition }\ (\ref{cstar})\ {\rm holds }.
\end{array}\right.
$$

\begin{remark}
\label{safarov} We note that if $\Omega $ is a domain of class $C^{0,1}$, {\it i.e.,} $\Omega$ is locally a subgraph of Lipschitz continuous functions, then inequality (\ref{cstar}) holds for any $\alpha > N/2$ (see {\it e.g.,} \cite{buda}, Netrusov and Safarov~\cite{netsa} and also Remark~\ref{ultrameyer} below). We also note that by the Min-Max Principle \cite[p.~5]{daheat} and by comparing with the Dirichlet Laplacian on a ball contained
in $\Omega $, condition (\ref{cstar}) does not hold for $\alpha \le N/2$ (no matter whether $\Omega $ is regular or not).
\end{remark}

In order to compare $L$  and $\tilde L$, we shall compare the respective pull-backs $H$ and $\tilde{H}$.
Since these act on different Hilbert spaces -- $L^2(\Omega , g\, dx)$ and $L^2(\Omega , \tilde{g}\, dx)$ -- we shall use the canonical unitary operator,
\[
w: L^2(\Omega , g\, dx) \longrightarrow L^2(\Omega , \tilde{g}\, dx) \;\; , \qquad u\mapsto wu\, ,
\]
defined as the multiplication by the function  $w:=g^{1/2}\tilde{g}^{-1/2}$.
We also introduce the multiplication operator $S$ on $(L^2(\Omega ))^N$  by the symmetric matrix
\begin{equation}
\label{esse}
 w^{-2}a^{-1/2}\tilde a a ^{-1/2}\, ,
\end{equation}
where the matrix $a$ is defined by (\ref{matrixa}) and the matrix $\tilde a$ is defined in the same way with $\tilde \phi$ replacing $\phi$.
If there is no ambiguity we shall denote the matrix (\ref{esse}) also by $S$.

As it will be clear in the sequel, in order to compare $H$ and $\tilde H$ we shall need an auxiliary operator.
Namely we shall consider the operator $T^*ST$, which is the non-negative self-adjoint operator in $L^2(\Omega,g\, dx)$ canonically associated with the sesquilinear form
\[
\int_{\Omega }( \tilde a \nabla u_1\cdot \nabla \bar u_2 )\tilde g    dx \; , \;\;\;  u_1,u_2\in {\mathcal{V}}\, .
\]
It is easily seen  that the operator $T^*ST$ is the pull-back to $\Omega$ via $\tilde\phi$ of the operator
\begin{equation}
\hat{L}:=\frac{\tilde g\circ \tilde\phi^{(-1)}}{g\circ \tilde\phi^{(-1)}}\tilde L
\label{lhat}
\end{equation}
defined on $L^2(\tilde\phi(\Omega))$. Thus in the sequel we shall deal with the operators $L,\tilde L$ and $\hat L$ and
the respective pull-backs $H$, $\tilde H$ and $T^*ST$. We shall repeatedly use the fact that these are pairwise unitarily equivalent.

We  denote by $\lambda_n[E]$, $n\in\N$, the eigenvalues of a non-negative self-adjoint operator $E$ with compact resolvent, arranged in non-decreasing order and repeated according to multiplicity, and by $\psi_n[E]$,
$n\in\N$, a corresponding orthonormal sequence of eigenfunctions.

\begin{lemma} Let {\rm (A)} be satisfied. Then
the operators $L$, $\tilde L$, $\hat L$, $H$, $\tilde H$ and $T^*ST$ have compact resolvents and the corresponding non-zero eigenvalues satisfy
the inequality
\begin{equation}
\label{elpeq}
\sum_{\lambda_n[E]\ne 0}\lambda_n[E]^{-\alpha } \le c c^{*},
\end{equation}
for $E=L,\tilde L, \hat L, H, \tilde H, T^*ST$,
where $c$ depends only on $N, \tau ,\theta$.
\label{elp}
\end{lemma}
{\em Proof.} We prove inequality (\ref{elpeq}) only for $E=T^*ST$, the other cases being similar.
Observe that the Rayleigh quotient corresponding to $T^*ST$ is given by
$$
\frac{\inprod{T^*ST u}{u}_{g}}{\inprod {u}{u}_g}=
\frac{\inprod{ST u}{ T u}_{g}}{\inprod {u}{u}_g}=\frac{\int_{\Omega }( \tilde a \nabla u\cdot \nabla \bar u )\tilde g    dx }{\int_{\Omega } |u|^2gdx },\  \ \ u\in {\mathcal{V}}.
$$
Then inequality (\ref{elpeq})  easily follows by
observing that
$$
 \tilde a \nabla u\cdot \nabla\bar u  \geq \theta^{-1} | (\nabla \tilde \phi )^{-1}\nabla u |^2\geq \theta^{-1}\tau ^{-2}|\nabla u|^2,
$$
\begin{equation}
| {\mathrm{det}} \nabla \phi | \leq
N! |\nabla\phi| ^{N} \label{in1}
\end{equation}
and using  the Min-Max Principle \cite[p.\,5]{daheat}.\hfill $\Box$

\

\section{Stability estimates for the resolvent and the eigenvalues}
\label{resolv}

The following lemma is based on the well-known commutation formula \eqref{deift} (see Deift~\cite{D}). By $\sigma (E)$ we denote the spectrum of an operator $E$.
\begin{lemma}\label{identity} Let {\rm (A)} be satisfied. Then
for all $\xi\in\C \setminus (\sigma (H)\cup \sigma (\tilde H)\cup \sigma(T^*ST))$,
\begin{equation}
(w^{-1}\tilde{H}w -\xi)^{-1} -(H-\xi )^{-1} =A_1+A_2+A_3+B\, ,
\label{identity1}
\end{equation}
 where
\begin{eqnarray*}
A_1&=&(1-w)(wT^*STw-\xi )^{-1}, \\
A_2&=&w(wT^*STw-\xi )^{-1}(1-w), \\
A_3&=& -\xi (T^*ST-\xi )^{-1}(w -w^{-1})(wT^*STw-\xi )^{-1}w, \\
B&=& T^*S^{1/2}(S^{1/2}TT^*S^{1/2}-\xi )^{-1}S^{1/2}(S^{-1}-I)(TT^*-\xi )^{-1}T\, .
\end{eqnarray*}
\end{lemma}

{\em Proof.} It suffices to prove (\ref{identity1}) for $\xi\neq 0$, since the case $\xi=0\not\in \sigma(H)\cup\sigma(\tilde H)\cup\sigma(T^*ST)$ is then obtained by letting $\xi\to 0$.

Recall that $T^*T=H$. Similarly  $\tilde{T}^{  \tilde{*}   }\tilde{T}=\tilde{H}$, where we have emphasized the dependence
of the adjoint operation on the specific inner-product used.
In this respect we note that the two adjoints of an operator $E$ are related by the conjugation relation $E^{\tilde {*}  }=w^2E^*w^{-2}$. This will allow us to use only $*$ and not $\tilde {*}$.

Observe that
\begin{equation}
\label{eich}
\tilde H= (\tilde a ^{1/2}\nabla ) ^{ \tilde {*}  } \tilde a ^{1/2}\nabla = w^2 (\tilde a ^{1/2}\nabla ) ^{*} w^{-2}\tilde a ^{1/2}\nabla   =w^2T^{*}ST\, .
\end{equation}
Therefore, by simple computations, we obtain
\begin{eqnarray*}
\lefteqn{(w^{-1}\tilde{H}w-\xi )^{-1}-(H-\xi)^{-1}} \\
\quad &=&w^{-1}(\tilde{H}-\xi)^{-1}w -(H-\xi)^{-1} \\
\quad &=& w^{-1}(w^2T^*ST-\xi)^{-1}w -(T^*T-\xi)^{-1} \\
\quad &=&w^{-1}(w^2T^*ST-\xi )^{-1}w -(T^*ST-\xi w^{-2})^{-1}w^{-1}\\
\quad &&+(T^*ST-\xi w^{-2})^{-1}w^{-1} -(T^*ST-\xi w^{-2})^{-1} \\
\quad &&+(T^*ST-\xi w^{-2})^{-1}-(T^*ST-\xi)^{-1} \\
\quad &&+ (T^*ST-\xi)^{-1} -(T^*T-\xi)^{-1} \\
\quad &=& A_1+A_2+A_3+((T^*ST-\xi)^{-1} -(T^*T-\xi )^{-1})\, .
\end{eqnarray*}
To compute the last term we  use the commutation formula
\begin{equation}
-\xi(E^*E-\xi)^{-1}+E^*(EE^*-\xi)^{-1}E =I
\label{deift}
\end{equation}
which holds for any closed and densely defined operator $E$, see \cite{D}. We write (\ref{deift}) first for $E=T$, then for $E=S^{1/2}T$, and then we subtract the two relations. After some simple calculations we obtain
$(T^*ST-\xi)^{-1} -(T^*T-\xi )^{-1}=B$, as required. \hfill $ \Box$\\

We now introduce a regularity property which will be important for our estimates. Sufficient conditions for its validity will be given in Section \ref{reg}.
\begin{definition}\label{growth}
Let $U$ be an open set in $\R^N$ and let $E$ be a non-negative self-adjoint operator on $L^2(U)$ with
compact resolvent and $\dom(E)\subset W^{1,2}(U)$. We say that $E$ satisfies property {\rm (P)} if
there exist $q_0>2, \gamma \geq 0 , C>0 $ such that the eigenfunctions $\psi_n[E]$  of $E$ satisfy the following conditions:\vspace{12pt}\\
\vspace{0cm}\hspace{2cm} $ \| \psi_n[E]\|_{L^{q_0}(U)}\leq C\lambda_n[E]^{\gamma}$  \hfill {\rm (P1)}\vspace{12pt}\\
and  \vspace{12pt}\\
\vspace{0cm}\hspace{2cm} $\| \nabla\psi_n[E]\|_{L^{q_0}(U)}  \leq C\lambda_n[E]^{\gamma +\frac{1}{2}  }$   \hfill {\rm (P2)}\vspace{12pt}\\
for all $n\in \N$ such that $\lambda_n[E]\ne 0$.
\end{definition}

\begin{remark}
\label{growthrem}
It is known that if $\Omega$, $A_{i,j}$ and $\Gamma$ are sufficiently smooth then for the operator $L$ in ({\ref{intro1}}),
subject to the boundary conditions ({\ref{intro2}}), property {\rm (P)} is satisfied with $q_0=\infty $ and $\gamma =N/4$; see Theorem~\ref{troy} and the proof of Theorem~\ref{thm:C11}.
\end{remark}

By interpolation it follows that if conditions {\rm (P1)} and {\rm (P2)} are satisfied then
\begin{equation}
\| \psi_n[E]\|_{L^q(U)}\leq C^{\frac{q_0(q-2)}{q(q_0-2)}}\lambda_n[E]^{\frac{q_0(q-2)\gamma}{q(q_0-2)}}\, , \;\;  \| \nabla\psi_n[E]\|_{L^q(U)}  \leq C^{\frac{q_0(q-2)}{q(q_0-2)}}\lambda_n[E]^{\frac{1}{2}+ \frac{q_0(q-2)\gamma}{q(q_0-2)}}\, ,
\label{regq}
\end{equation}
for all $q\in [2,q_0]$.

\

In the sequel we  require that property (P) is satisfied by the operators $H$, $\tilde H$ and $T^*ST$ which, according to the following lemma, is equivalent to requiring   that property (P) is satisfied by the operators $L$, $\tilde L$ and $\hat L$ respectively.

\begin{lemma}
\label{equivp}
Let {\rm (A)} be satisfied.  Then the operators   $H$, $\tilde H$ and $T^*ST$ respectively, satisfy property
{\rm (P)} for some $q_0>2$ and $\gamma \geq 0$ if and only if the operators $L$, $\tilde L$ and $\hat L$ respectively, satisfy property
{\rm (P)} for the same $q_0$ and $\gamma$.
\end{lemma}

Let $E$ be a non-negative self-adjoint operator on a Hilbert space the spectrum of which consists of isolated positive eigenvalues of finite multiplicity and may also contain zero as an eigenvalue of possibly infinite multiplicity.
Let $s>0$. Given a function $g:\sigma(E)\to\C$ we define
\begin{eqnarray*}
&& |g(E)|_{p,s}= \Big(\sum_{\lambda_n[E]\neq 0} |g(\lambda_n[E])|^p \lambda_n[E]^s\Big)^{1/p}, \; 1\leq p<\infty\, ,\\
&& |g(E)|_{\infty,s}= \sup_{\lambda_n[E]\neq 0} |g(\lambda_n[E])|\, ,
\end{eqnarray*}
where, as usual, each positive eigenvalue is repeated according to its multiplicity.

The next lemma involves the Schatten norms $\|\cdot\|_{\cC^r}$, $1\leq r\leq\infty$. For a compact operator $E$ on a Hilbert space they are defined by $\|E\|_{\cC^r}=(\sum_n \mu_n(E)^r)^{1/r}$, if $r<\infty$, and $\|E\|_{\cC^{\infty}}=\|E\|$,
where $\mu_n(E)$ are the singular values of $E$, {\it i.e.,} the non-zero eigenvalues of $(E^*E)^{1/2}$;
recall that the Schatten space $\cC^r$, defined as the space of those compact operators for which the Schatten norm $\|\cdot\|_{\cC^r}$ is finite, is a Banach space; see Reed and Simon~\cite{RS} or Simon~\cite{S} for details.

Let $F:=TT^*$. Recall that $\sigma(F)\setminus\{0\}=\sigma(H)\setminus\{0\}$, see \cite{D}. In the next lemma, $g(H)$ and $g(F)$ are operators defined in the standard way by functional calculus. The following lemma is a variant of Lemma 8 of \cite{bar3}.
\begin{lemma}
\label{interp} Let $q_0>2$, $\gamma\geq 0$, $p\geq q_0/(q_0-2)$, $2\leq r <\infty$ and $s=2q_0\gamma/[p(q_0-2)]$.
Then the following statements hold:
\begin{itemize}
\item[(i)] If the eigenfunctions of $H$ satisfy {\rm (P1)} then
for any measurable function $R:\Omega\to\C$ and any function $g:\sigma(H)\to\C$ we have
\begin{equation}
\| Rg(H)\|_{\cC^{r}}\leq   \| R\|_{L^{pr}(\Omega)}\Big( |\Omega |^{-\frac{1}{pr}} |g(0)| +  C^{\frac{2q_0}{pr(q_0-2)}} |g(H)|_{r,s}\Big)\, .
\label{mainrembis}
\end{equation}
\item[(ii)] If the eigenfunctions of $H$ satisfy {\rm (P2)} then
for any measurable matrix-valued function $R$ on $\Omega$ and any function $g:\sigma(F)\to\C$ such that if $0\in \sigma (F)$ then $g(0)=0$,  we have
\begin{equation}
\| Rg(F)\|_{\cC^{r}}\leq   C^{\frac{2q_0}{pr(q_0-2)}} \| a \|_{L^{\infty}(\Omega ) }^{\frac{1}{r}}\,  \|R\|_{L^{pr}(\Omega)}|g(F)|_{r,s}\, .
\label{1}
\end{equation}
\end{itemize}
\end{lemma}
{\em Proof.} We  only prove statement {\it (ii)} since the proof of {\it (i)} is simpler.
It is enough to consider the case where $R$ is bounded and $g$ has finite support: the general case will then follow by approximating $R$ in $\|\cdot \|_{L^{pr}(\Omega )}$ by a sequence $R_n,\, n\in\N$, of bounded matrix-valued  functions and $g$ in $|\cdot |_{r,s}$ by a sequence $g_n$, $n\in \N$,  of functions with finite support, and observing that by (\ref{1}) the sequence $R_ng_n(F), \, n\in\N$, is then a Cauchy sequence in $\cC^r$.

Since $R$ is bounded and $g$ has finite support  $Rg(F)$ is compact, hence inequality (\ref{1}) is trivial for $r=\infty $.
Thus it is enough to prove (\ref{1}) for $r=2$
since  the general case will then follow by interpolation (cf. \cite{S}).
It is easily seen that $z_n:=T\psi_n[H]/ \|T\psi_n[H]\| = \lambda_n[H]^{-1/2}T\psi_n[H] $, for all $n\in \N$ such that $\lambda_n[H]\ne 0$, are orthonormal eigenfunctions of $F$, $Fz_n=\lambda_n[H]z_n$, $n\in\N$, and
${\rm span}\{z_n\}=\kernel (F)^{\perp}$. Hence, since $g(0)=0$,
\begin{eqnarray}
\|R g(F)\|_{\cC^2}^2&=&\sum_{n=1}^{\infty}\|Rg(F)z_n\|_{L^2(\Omega )}^2 \label{hilbert}\\
&=&\sum_{n=1}^{\infty}|g(\lambda_n[H])|^2\|Rz_n\|_{L^2(\Omega )}^2 \nonumber \\
&=&\sum_{n=1}^{\infty}\lambda_n[H]^{-1}|g(\lambda_n[H])|^2\|Ra^{1/2}\nabla \psi_n[H]\|_{L^2(\Omega )}^2\nonumber\\
&\leq& \|a^{1/2}\|_{L^{\infty}(\Omega )}^2 \|R\|_{L^{2p}(\Omega )}^2\sum_{n=1}^{\infty}\lambda_n[H]^{-1}|g(\lambda_n[H])|^2\|\nabla \psi_n[H]\|^2_{L^{2p/(p-1)}(\Omega )}\nonumber\\
&\leq&C^{\frac{2q_0}{p(q_0-2)}}\|a^{1/2}\|_{L^\infty (\Omega )}^2\|R\|_{L^{2p}(\Omega )}^2\sum_{n=1}^{\infty}|g(\lambda_n[H])|^2\lambda_n[H]^{\frac{2q_0\gamma}{p(q_0-2)}},\nonumber
\end{eqnarray}
where for the last inequality we have used (\ref{regq}).
This proves (\ref{1}) for $r=2$, thus completing the proof of the lemma. $\hfill \Box$\\

Recall that $\delta _p(\phi , \tilde \phi )$,  $1\le p\le\infty $,
is defined in (\ref{deltaphi}).

\begin{theorem} {\bf (stability of resolvents)} Let {\rm (A)} be satisfied.
Let $\xi\in \C\setminus (\sigma (H)\cup \{0\} )$.
Then the following statements hold:
\begin{itemize}
\item[(i)] There exists $c_1>0$ depending only on $N,  \tau ,\theta , \alpha , c^*$ and $\xi $ such that if
$
\delta _{\infty}(\phi ,\tilde \phi )$ $ \leq c_1^{-1},
$
then $\xi \notin \sigma (\tilde H)$ and
\begin{equation}
\|(w^{-1}\tilde{H}w -\xi)^{-1} -(H-\xi )^{-1}\|_{\cC^{\alpha }}\leq c_1  \delta _{\infty}(\phi ,\tilde \phi ).
\label{res1}
\end{equation}
\item[(ii)]
Let in addition  {\rm (P)} be satisfied by the operators $H$, $\tilde H$ and $T^*ST$ for the same $q_0$, $\gamma $ and $C$. Let $p\geq q_0/(q_0-2)$ and  $r \geq  \max \{2,    \alpha  +  \frac{2q_0\gamma}{p(q_0-2)}  \}$. Then there exists  $c_2>0$ depending only on $N , \tau ,\theta , \alpha, c^*, r, p,q_0,C, \gamma ,  |\Omega |$ and $\xi $ such that if
$
\delta _{pr}(\phi ,\tilde \phi )   \leq c_2^{-1},
$
then $\xi \notin \sigma (\tilde H)$ and
\begin{equation}
\|(w^{-1}\tilde{H}w -\xi)^{-1} -(H-\xi )^{-1}\|_{\cC^r}\leq c_2
\delta _{pr}(\phi ,\tilde \phi ).
\label{res2}
\end{equation}
\end{itemize}
\label{mainthm}
\end{theorem}

\begin{remark}
Let
$
s=[q_0/(q_0-2)] \max\{ 2, \alpha+2\gamma\}.
$
It follows by Theorem \ref{mainthm} (ii) (choosing $p=q_0/(q_0-2)$) that
if
$
 \delta_s(\phi ,\tilde \phi) \leq c_2^{-1},
$
then $\xi \notin \sigma (\tilde H)$ and
\begin{equation}
\|(w^{-1}\tilde{H}w -\xi)^{-1} -(H-\xi )^{-1}\|\leq c_2
  \delta_s(\phi ,\tilde \phi).
\label{res21}
\end{equation}
\end{remark}

\begin{remark}
\label{bestrange}
As we shall see in Section 7, the best range for $s$ in (\ref{res21}) used in our applications is $s >N$; this corresponds to the case in which
$q_0=\infty $, $\gamma =N/4$ and $\alpha >N/2$. See Remarks~\ref{safarov}, \ref{growthrem}.
\end{remark}

{\em Proof of Theorem \ref{mainthm}.}
In this proof we denote by $c$
a positive constant depending only on $N, \tau ,\theta, \alpha$ and $c^*$ the value of which  may change along the proof; when dealing with statement $(ii)$ constant $c$ may depend also on $r, p, q_0, C, \gamma , |\Omega |$.
We divide the proof into two steps.

{\em Step 1.} We assume first that $\xi\not\in \sigma (\tilde H)\cup \sigma(T^*ST)$ and we set
\[
d_{\sigma}(\xi)=\dist(\xi, \sigma(H)\cup\sigma(\tilde{H})\cup\sigma(T^*ST)).
\]
In this first step we shall prove (\ref{res1}) and (\ref{res2})  without any smallness assumptions on $\delta_{\infty }(\phi ,\tilde \phi )$,
 $\delta_{pr }(\phi ,\tilde \phi )$ respectively.

We first prove (\ref{res1}). We shall use Lemma~\ref{identity} and to do so we first estimate the terms $A_1, A_2, A_3$ in identity (\ref{identity1}).
Clearly we have that
\begin{equation}
\label{crtildeh1}
\frac{\lambda_n [\tilde H] }{| \lambda_n[\tilde H] -\xi |}\le \left(1+\frac{|\xi |}{d(\xi , \sigma (\tilde H))}  \right)\, .
\end{equation}
Since the eigenvalues of the operator $wT^*STw$ coincide with the eigenvalues of  $\tilde H$ (see (\ref{eich})), it follows that
\begin{eqnarray}
\label{crtildeh2}
\|(wT^*STw-\xi )^{-1}\|_{\cC^{\alpha }}^{\alpha }
& =& \sum_{n=1}^{\infty }\frac{1}{|\lambda_n[\tilde H]-\xi |^{\alpha }}
\nonumber \\
&
\le & \frac{1}{|\xi |^{\alpha }}+    \left( 1+\frac{|\xi |}{ d(\xi , \sigma (\tilde H))} \right)^{\alpha } \sum_{\lambda_n[\tilde H]\ne 0}\lambda_n[\tilde H]^{-\alpha }\nonumber \\
& =&  \frac{1}{|\xi |^{\alpha }}+ c    \left( 1+\frac{|\xi |}{ d(\xi , \sigma (\tilde H))} \right)^{\alpha } .
\end{eqnarray}
By (\ref{in1}) and by observing that
\begin{equation}
\big|{\mathrm{det}} \nabla \phi-{\mathrm{det}} \nabla\tilde{\phi}\big|
 \leq  N! N\,\big|\nabla\phi-\nabla\tilde{\phi}\big|\,
\max\left\{
\,\left|\nabla\phi\right|\,,
\,\big|\nabla\tilde{\phi}\big|\,
\right\}^{N-1}\label{in3}
\end{equation}
it follows that
\begin{equation}
|1-w  | , |w-w^{-1}  |  \le   c  |\nabla \phi -\nabla\tilde \phi  |.  \label{crtildeh3}
\end{equation}
Combining inequalities (\ref{crtildeh2}) and  (\ref{crtildeh3}) we obtain
that
\begin{eqnarray} & &
\| A_1\|_{\cC^{\alpha }}, \| A_2\|_{\cC^{\alpha }}\le c \left( 1+\frac{1}{|\xi |} +\frac{|\xi |}{d_{\sigma }(\xi )  }  \right)   \|\nabla \phi -\nabla\tilde \phi  \|_{L^{\infty }(\Omega )}\, ,   \nonumber  \\
& & \| A_3\|_{\cC^{\alpha }}\le c \left( \frac{1+|\xi |}{d_{\sigma}(\xi)} + \frac{|\xi |^2}{d_{\sigma }(\xi )^2}  \right)   \|\nabla \phi -\nabla\tilde \phi  \|_{ L^{\infty }(\Omega )}\, .
\label{a1a2a3}
\end{eqnarray}

We now estimate the term  $B$ in (\ref{identity1}). We recall that $F=TT^*$ and we set  $F_S=S^{1/2}TT^*S^{1/2}$. Then, by polar decomposition, there exist partial isometries
$Y,Y_S: L^2(\Omega,g\, dx)\to (L^2(\Omega,g\, dx))^N$ such that $T=F^{1/2}Y$ and $S^{1/2}T=F_S^{1/2}Y_S$.
We then have
$$
B=Y_S^* F_S^{1/2}(F_S-\xi)^{-1}S^{1/2}(S^{-1}-I)(F-\xi)^{-1}F^{1/2}Y\, .
$$
Hence, by H\"{o}lder's inequality for the Schatten norms (see \cite[p.~41]{RS}) it follows that
\begin{equation}\label{s-2}
\| B\|_{\cC^{\alpha }}\leq \| F_S^{1/2}(F_S-\xi)^{-1}\|_{\cC^{2{\alpha }}}
\|S^{1/2}(S^{-1}-I)\|_{ L^{\infty }(\Omega )} \|(F-\xi)^{-1}F^{1/2}\|_{\cC^{2{\alpha }}}\, .
\end{equation}
Since $\sigma(F)\setminus \{0\}=\sigma(H)\setminus \{0\}$, we may argue as before and obtain
\begin{eqnarray}
\label{s-1}
\|(F-\xi )^{-1}F^{1/2}\|_{\cC^{2{\alpha }}}^{2{\alpha }} & \leq& c   \left(1+\frac{|\xi |}{ d(\xi , \sigma (H)) } \right)^{2{\alpha }}, \nonumber \\
\|(F_S-\xi )^{-1}F_S^{1/2}\|_{\cC^{2{\alpha }}}^{2{\alpha }}&
\le & c \left(1+\frac{|\xi |}{ d(\xi , \sigma (T^*ST)) } \right)^{2\alpha }.
\end{eqnarray}
Now, one easily sees that
\begin{eqnarray}
\label{s1}
|S^{-1}-I| & \leq  &   |(w^2-1)a^{1/2}\tilde{a}^{-1}a^{1/2}|+
|a^{1/2}(\tilde{a}^{-1}-a^{-1})a^{1/2}|\nonumber \\ & \leq &  c (|\nabla\phi-\nabla\tilde{\phi} |+ |A\circ \phi -A\circ \tilde \phi|).
\end{eqnarray}
Combining (\ref{s-2}), (\ref{s-1}) and  (\ref{s1})  we conclude that
\begin{equation}
\label{b}
\|B\|_{\cC^{\alpha }}\leq c \left(1+\frac{|\xi |}{ d_{\sigma}(\xi) } \right)^2 \delta_{\infty }(\phi ,\tilde\phi).
\end{equation}
By Lemma~\ref{identity} and estimates (\ref{a1a2a3}) and (\ref{b}), we deduce that

\begin{equation}
\label{35bis}
\|(w^{-1}\tilde{H}w -\xi)^{-1} -(H-\xi )^{-1}\|_{\cC^{\alpha }}\leq c_1\left( 1+\frac{1}{|\xi|}+\frac{1}{d_{ \sigma   }(\xi )}+\frac{|\xi|^2}{d_{\sigma    }
(\xi)^2}\right) \delta _{\infty}(\phi ,\tilde \phi ).
\end{equation}

We now prove $(\ref{res2})$.
In order to estimate $A_1$, $A_2$, $A_3$ we use estimate (\ref{mainrembis}) and we get
\begin{eqnarray}
& &
\| A_1\|_{\cC^r}, \| A_2\|_{\cC^r} \le c
\left( 1+\frac{1}{|\xi|}+\frac{|\xi |}{d_{\sigma }(\xi )}\right)
\|\nabla\phi-\nabla\tilde{\phi}\|_{L^{pr}(\Omega )}\, ,\nonumber \label{a1a2a3abis}\\  & &
\| A_3\|_{\cC^r}  \le c
\left(   \frac{1+|\xi |}{d_{\sigma}(\xi)} + \frac{|\xi |^2}{d_{\sigma }(\xi )^2}  \right)
\|\nabla\phi-\nabla\tilde{\phi}\|_{ L^{pr}(\Omega )}\label{a1a2a3bis}\, .
\end{eqnarray}

We now estimate $B$.
We shall assume without loss of generality that $S^{-1}-I\geq 0$. Thus, in order to estimate the $\cC^r$ norm of $B$, we shall estimate the $\cC^{2r}$ norms of  $F_S^{1/2}(F_S-\xi)^{-1}S^{1/2}(S^{-1}-I)^{1/2}$ and $(S^{-1}-I)^{1/2}(F-\xi)^{-1}F^{1/2}$.
By Lemma \ref{interp} it follows that
\begin{eqnarray}
\lefteqn{\|(S^{-1}-I)^{1/2}(F-\xi)^{-1}F^{1/2}\| _{\cC^{2r}}^{2r}} \nonumber \\
& & \quad \le  c \| (S^{-1}-I)^{1/2}  \|_{L^{2pr}(\Omega ) }^{2r} \sum_{n=1}^{\infty }\biggl| \frac{\lambda_n[H]}{\lambda_n[H] -\xi }  \biggr|^{2r}\lambda_n[H]^{  \frac{2q_0\gamma}{p(q_0-2)}   -r}\nonumber\\
& &\quad \le  \quad c   \| S^{-1}-I  \|_{L^{pr}(\Omega )}^{r}
 \left(1+\frac{|\xi |}{ d_{\sigma }(\xi ) } \right)^{2r} .
\end{eqnarray}
The same estimate  holds also for the operator $F_S^{1/2}(F_S-\xi)^{-1}S^{1/2}(S^{-1}-I)^{1/2}$. Thus by H\"{o}lder inequality for the Schatten norms it follows that
\begin{equation}
\label{bbis}
\| B\|_{\cC^r}\le c \left(1+\frac{|\xi |}{ d_{\sigma }(\xi ) } \right)^{2}\| S^{-1}-I  \|_{L^{pr}(\Omega )}.
\end{equation}
By Lemma~\ref{identity} and combining estimates and (\ref{s1}), (\ref{a1a2a3bis}) and (\ref{bbis})  we deduce that

\begin{equation}
\label{38bis}
\|(w^{-1}\tilde{H}w -\xi)^{-1} -(H-\xi )^{-1}\|_{\cC^r}\leq c_1\left( 1+\frac{1}{|\xi|}+\frac{1}{d_{ \sigma   }(\xi )}+\frac{|\xi|^2}{d_{\sigma    }
(\xi)^2}\right) \delta _{pr}(\phi ,\tilde \phi ).
\end{equation}

{\em Step 2.} We prove statement {\it (i)}. First of all we  prove that there exists $c>0$ such that if
\begin{equation}
\label{claim}
\delta_{\infty }(\phi ,\tilde\phi)   < \frac{d(\xi , \sigma (H) ) }{c(1+|\xi |^2+d(\xi , \sigma (H))^2)  }
\end{equation}
then  $\xi \notin \sigma (\tilde H) \cup \sigma (T^*ST)$ and
\begin{equation}
d(\xi ,  \sigma (\tilde H)  ), \
d(\xi , \sigma (T^*ST))>\frac{d(\xi , \sigma (H))}{2}.
\label{claim05}
\end{equation}

We begin with $T^*ST$. By recalling that $B=(T^*ST-\xi)^{-1}-(T^*T-\xi)^{-1}   $ (see the proof of Lemma~\ref{identity}), by estimate (\ref{b}) with $\xi =-1$ and by inequality (\ref{abstsimon}) it follows that there exists $C_1>0$ such that for all $n\in\N$
\begin{equation}
\label{claim1}
\left|  \frac{1}{\lambda_n [T^*ST ]+1}- \frac{1}{\lambda_n[H]+1}  \right|\le C_1 \delta_{\infty }(\phi ,\tilde\phi).
\end{equation}

Assume that $n\in \N$ is such that $\lambda_n [T^*ST ]\le |\xi |+d(\xi , \sigma (H))$; by (\ref{claim1}) it follows that if
$$
C_1(1+|\xi |+d(\xi , \sigma (H)))\delta_{\infty }(\phi ,\tilde\phi)  <
\frac{ |\xi |+d(\xi , \sigma (H))}{2(|\xi |+d (\xi, \sigma (H))  )+1}
$$ then
\begin{eqnarray}
\lambda _n[H]  &  \le & \frac{|\xi |+d (\xi , \sigma (H))  +C_1[1+|\xi | +d(\xi , \sigma (H))  ] \delta_{\infty }(\phi ,\tilde\phi) }{1-C_1[1+|\xi | +d (\xi , \sigma (H))] \delta_{\infty }(\phi ,\tilde\phi)  } \nonumber\\ & \le & 2(|\xi |+d(\xi , \sigma (H))),
\label{claim2}
\end{eqnarray}
(the elementary inequality $(A+t )(1-t)^{-1}<2A$ if $0<t< A(2A+1)^{-1}$ was used).
Thus by (\ref{claim1}) and (\ref{claim2}) it follows that if
$$\delta_{\infty }(\phi,\tilde\phi) \le
\frac{d(\xi , \sigma (H))}{    2C_1[ 1+|\xi | +d(\xi , \sigma (H))]  [1+2(|\xi | +d(\xi , \sigma (H)))  ]} $$
then
\begin{eqnarray}\lefteqn{
 |\xi -\lambda_n[T^*ST]|   } \nonumber  \\ &  \geq & |\xi -\lambda_n[H]|-|\lambda_n[H]-\lambda_n[T^*ST]| \nonumber \\
 & \geq &  d(\xi, \sigma (H))-C_1[1+|\xi |+d(\xi, \sigma (H))]
[1+2(|\xi |+d(\xi, \sigma (H)))]
\delta_{\infty }(\phi ,\tilde\phi)
\nonumber \\ \label{claim3}
& \geq & \frac{d(\xi , \sigma (H))}{2}\, ,
\end{eqnarray}
for all $n\in \N$ such that $\lambda_n [T^*ST ]\le |\xi |+d(\xi , \sigma (H))$. Thus inequality (\ref{claim05}) for $d(\xi , \sigma (T^*ST))$  follows by (\ref{claim3}) and by observing that
if $n\in \N$ is such that $\lambda_n [T^*ST ]> |\xi |+d(\xi , \sigma (H))$ then $|\xi -\lambda_n[T^*ST]| >d(\xi,\sigma(H))$.
Inequality (\ref{claim05}) for $d(\xi , \sigma (\tilde H))$ can be proved in the same way. Indeed it suffices to  observe that by {\it Step 1} there exists $ C_2>0$ such that
\begin{equation}
\label{claim1'}
\left|  \frac{1}{\lambda_n [\tilde H]+1}- \frac{1}{\lambda_n[H]+1}  \right|\le  C_2 \delta_{\infty }(\phi ,\tilde\phi);
\end{equation}
we then proceed exactly as above.

By (\ref{35bis}) and (\ref{claim05}) it follows that there exists $c>0$ such that if
\begin{equation}
\delta _{\infty}(\phi ,\tilde \phi ) \leq \frac{d_{   }(\xi , \sigma (H) )}{c (1+|\xi|^2+d_{  }(\xi , \sigma (H) )^2)},
\label{res1a}
\end{equation}
then $\xi \notin \sigma (\tilde H)\cup \sigma (T^*ST)$ and
\begin{equation}
\label{35bisbis}
\|(w^{-1}\tilde{H}w -\xi)^{-1} -(H-\xi )^{-1}\|_{\cC^{\alpha }}\leq c \left( \! 1+\frac{1}{|\xi|}+\frac{1}{d(\xi
, \sigma (H) )}+\frac{|\xi|^2}{d
(\xi  , \sigma (H) )^2}\right) \delta _{\infty}(\phi ,\tilde \phi ) .
\end{equation}
This completes the proof of statement $(i)$.

The argument above  works word by word also for the proof of statement $(ii)$, provided  that $\delta_{\infty }(\phi ,\tilde\phi)$ is replaced
by $ \delta_{pr } (\phi ,\tilde\phi)   $. \hfill $\Box$

\

\begin{remark}
\label{rem:distance}
The proof of Theorem \ref{mainthm} gives  some information about the dependence of the  constants $c_1, c_2$ on $\xi$ which will be useful in the sequel. For instance, in the case of statement (i), in fact it was proved that there exists $c$ depending only on $N, \tau ,\theta , \alpha , c^*$ such that
if ({\ref{res1a}}) holds then (\ref{35bisbis}) holds.
Exactly the same holds for statement (ii) where $c$ depends also on $r, p,q_0,C, \gamma , |\Omega |$.
Moreover, for such $\phi , \tilde\phi$,  if $0\not\in \sigma(H)$ then $0\not\in\sigma(\tilde H)$ and the summand $1/|\xi| +1/d_{}(\xi , \sigma (\Omega ))$ can be removed from the right-hand side of  (\ref{35bisbis}); furthermore, in this case statements (i) and (ii) also hold for $\xi=0$. This can be easily seen by looking closely at the proofs of (\ref{a1a2a3}) and (\ref{a1a2a3abis}).
\end{remark}

\begin{remark} By the proof of Theorem \ref{mainthm} it also follows that
for a fixed $\xi \in \C \setminus [0, \infty [$  no smallness conditions on $\delta_{\infty } (\phi ,\tilde\phi)$, $\delta_{pr } (\phi ,\tilde\phi)$ are required for the
validity of (\ref{res1}), (\ref{res2}) respectively.
\end{remark}

\begin{theorem} {\bf (stability of eigenvalues)} Let   {\rm (A)} be satisfied. Then
the following statements hold:
\begin{itemize}
\item[(i)] There exists $c_1>0$ depending only on $N, \tau ,\theta , \alpha$ and $c^*$ such that if
$
\delta_{\infty }(\phi ,\tilde\phi ) \leq c_1^{-1},
$
then
\begin{equation}
\left( \sum_{n=1}^{\infty}\left| \frac{1}{
{\lambda}_n[\tilde L] +1} - \frac{1}{{\lambda}_n[L ] +1}\right|^{\alpha }\right)^{1/{\alpha }} \leq c_1  \delta_{\infty }(\phi ,\tilde\phi ).
\label{res1'}
\end{equation}
\item[(ii)]
Let in addition {\rm (P)} be satisfied by the operators $L$, $\tilde L$ and $\hat L$ for the same $q_0>2$, $\gamma \geq 0 $ and $C>0$. Let $p\geq q_0/(q_0-2)$ and  $r\geq  \max \{2,    \alpha  + \frac{2q_0\gamma}{p(q_0-2)}   \}$. Then there exists  $c_2>0$ depending only on $N,\tau ,\theta , \alpha, c^*, r, p,q_0,C, \gamma $ and $|\Omega |$  such that if
$
\delta_{pr }(\phi ,\tilde\phi ) \le c_2^{-1},
$
then
\begin{equation}
\left( \sum_{n=1}^{\infty}\left| \frac{1}{
{\lambda}_n[\tilde L] +1} - \frac{1}{{\lambda}_n[L] +1}\right|^r\right)^{1/r}
\leq c_2 \delta_{pr }(\phi ,\tilde\phi ).
\label{res2'}
\end{equation}
\end{itemize}
\label{thm:series}
\end{theorem}
{\em Proof.} The theorem follows by Theorem \ref{mainthm} and  by applying the  inequality
\begin{equation}
\label{abstsimon}
\left( \sum_{n=1}^{\infty}\left| \frac{1}{
{\lambda}_n[E_1] +1} - \frac{1}{{\lambda}_n[E_2] +1}\right|^r\right)^{1/r}\leq
\|( E_1 +1)^{-1} -(E_2+1 )^{-1}\|_{\cC^r},
\end{equation}
with $E_1= w^{-1}\tilde H w$, $E_2=H$ (see \cite[p. 20]{S}).
$\hfill \Box$

\begin{remark}
We note that in the case of Dirichlet boundary conditions, {\it i.e., } ${\mathcal{V}}=W^{1,2}_0(\Omega )$, inequality (\ref{res1'}) directly follows by
 \cite[Lemma~6.1]{bula} the proof of which is based on the Min-Max Principle.
\end{remark}

\section{Stability estimates for eigenfunctions}
\label{eigenprsec}

\begin{definition}
Let $E$ be a non-negative self-adoint operator with compact resolvent on a Hilbert space ${\mathcal{H}}$. Given a finite subset $G$ of $\N$ we denote by $P_G(E)$ the orthogonal projector from ${\mathcal{H}}$ onto
the linear space generated by all the eigenfunctions corresponding to the eigenvalues $\lambda_k[E]$ with
$k\in G$.
\end{definition}

Observe that the dimension of the range of $P_G(E)$ coincides with the number of elements of $ G$ if and only if no eigenvalue with index in $G$ coincides with an eigenvalue with index in $\N \setminus G$; this will always be the case in what follows.

In the following statements it is understood that whenever $n=1$ the term $\lambda_{n-1}$ has to be dropped.

\begin{theorem} Let {\rm (A)} be satisfied.
Let $\lambda$ be a non-zero eigenvalue of $H$ of multiplicity $m$, let $n\in\N$ be such that $\lambda=\lambda_n[H]
=\ldots =\lambda_{n+m-1}[H]$, and let $G=\{n,n+1,\ldots,n+m-1\}$.  Then the following statements hold:

\begin{itemize}

\item[(i)] There exists $c_1>0$ depending only on $N,\tau ,\theta ,\alpha ,  c^*, \lambda_{n-1}[H], \lambda , \lambda _{n+m}[H] $ such that if $ \delta_{\infty }(\phi ,\tilde\phi )  \leq c_1^{-1},$
then $\dim \ran P_G(w^{-1}\tilde Hw)=m$ and
\begin{equation}
\|P_G(H)-P_G(w^{-1}\tilde Hw)\| \leq c_1  \delta_{\infty  }(\phi ,\tilde\phi ).
\label{proj1}
\end{equation}
\item[(ii)]
Let in addition {\rm (P)} be satisfied by the operators $H$, $\tilde H$ and $T^*ST$ for the same $q_0$, $\gamma$ and $C$. Let  $s=[q_0/(q_0-2)] \max\{ 2, \alpha+2\gamma\}$.
Then there exists  $c_2>0$ depending only on $N,\tau, \theta ,\alpha, c^*, q_0, C, \gamma ,  |\Omega |$ $   \lambda_{n-1}[H], \lambda , $ $ \lambda _{n+m}[H] $ such that if
$ \delta_{s }(\phi ,\tilde\phi )\leq c_2^{-1}$
then $\dim \ran P_G(w^{-1}\tilde Hw)=m$ and
\begin{equation}
\|P_G(H)-P_G(w^{-1}\tilde Hw)\| \leq c_2 \delta_{s }(\phi ,\tilde\phi ).
\label{proj2}
\end{equation}
\end{itemize}
\label{thm:projectors}
\end{theorem}
{\em Proof.} We set $\rho =\frac{1}{2}\dist(\lambda,(\sigma(H)\cup\{0\})\setminus
\{\lambda\})$ and $\lambda^*=\lambda $ if $\lambda$ is the first non-zero eigenvalue of $H$, and $\lambda^*=\lambda_{n-1}[H]$ otherwise.

By  Theorem \ref{thm:series} $(i)$ it follows that
\begin{equation}
|\lambda_k[H]-\lambda _k[\tilde H]| \leq c(\lambda_k[H]+1)(\lambda_k[\tilde H]+1)
\delta_{\infty  }(\phi ,\tilde\phi )\, .
\label{proj3}
\end{equation}
This implies that there exists $c>0$ such that if $\delta_{\infty  }(\phi ,\tilde\phi )< c^{-1}\lambda_k[H]/(\lambda_k[H]+1)^2$, then $\lambda_k[\tilde H]\leq
2\lambda_k[H]$. This together with (\ref{proj3}) implies the existence of $c>0$ such that if
$\delta_{\infty  }(\phi ,\tilde\phi )< c^{-1}\min\{\rho ,\lambda_k[H]\}/(\lambda_k[H]+1)^2$ then $|\lambda_k[H]-\lambda_k[\tilde H]|<\rho /2$. Applying this inequality for $k=n-1, \dots , n+m$, we deduce that if
$$\delta_{\infty  }(\phi ,\tilde\phi ) <\frac{\min\{\rho ,\lambda^*\}}{c(\lambda_{n+m}[H] +1)^2},$$ then
\begin{eqnarray}
& & |\lambda_k[\tilde H]-\lambda  |\leq \rho /2\, \ , \ \forall\ k\in G\, \nonumber \\
& & |\lambda_k[\tilde H]-\lambda  |\geq 3\rho /2 \, \ ,  \ \forall\ k\in \N\setminus G\, .
\label{proj4}
\end{eqnarray}
Hence  $\dim \ran P_G(w^{-1}\tilde Hw)=m$ and by  the well-known Riesz formula we have that
\begin{eqnarray}
& &
P_G[H]=-\frac{1}{2\pi i}\int_{\Gamma}(H-\xi )^{-1}d\xi\, , \\
& &  P_G[w^{-1}\tilde Hw]=-\frac{1}{2\pi i}\int_{\Gamma}(w^{-1}\tilde Hw-\xi )^{-1}d\xi\, ,
\end{eqnarray}
where $\Gamma (\theta )= \lambda +\rho  e^{i \theta }$, $0\le \theta < 2\pi $. Hence
\begin{equation}\label{rieszest} \| P_G[H] -P_G[w^{-1}\tilde Hw] \|\le \rho \sup_{\xi \in \Gamma} \|(w^{-1}\tilde{H}w -\xi)^{-1} -(H-\xi )^{-1}\|  .
\end{equation}
Let $c_1$ be as in Theorem \ref{mainthm} $(i)$. By Theorem \ref{mainthm} $(i)$ and Remark \ref{rem:distance}
and by observing that $\lambda -\rho \le |\xi |\le  \lambda +\rho  $ and  $1/|\xi |\le 1/ \rho $  for all  $\xi \in \Gamma$,
 it follows that if
$$
\delta_{\infty }(\phi ,\tilde\phi) < \frac{\rho}{c_1(1+\lambda^2_{n+m}[H] +\rho^2 )}
$$
then
\begin{equation}
\label{rieszest1}
\|(w^{-1}\tilde{H}w -\xi)^{-1} -(H-\xi )^{-1}\| \leq c_1\left( 1+\frac{1}{\rho}+\frac{\lambda^2}{\rho ^2}   \right)\delta _{\infty }(\phi ,\tilde\phi ).
\end{equation}

The proof of statement $(i)$ then follows by combining (\ref{rieszest}) and (\ref{rieszest1}).
The proof of statement $(ii)$ is similar. \hfill $\Box $

\

\begin{remark}
The proof of Theorem \ref{thm:projectors} gives some information about the dependence of the constants $c_1,c_2$ on
$ \lambda_{n-1}[H], \lambda , \lambda _{n+m}[H] $ which will be useful in the sequel. For instance in the case of statement (i) in fact we have proved that
there exists $c>0$ depending only on $N,\tau ,\theta, \alpha, c^*$ such that  if
$$
 \delta_{\infty }(\phi ,\tilde\phi )  \leq \frac{\min\{\rho ,\lambda^*\}}{c(1+\rho ^2+\lambda_{n+m}[H]^2) }
$$
then
\begin{equation}
\|P_G(H)-P_G(w^{-1}\tilde Hw)\| \leq c\left(1+\rho+ \frac{\lambda^2}{\rho }\right)  \delta_{\infty  }(\phi ,\tilde\phi ).
\end{equation}
Exactly the same is true for statement (ii) where $c$ depends also on $q_0, C, \gamma , |\Omega |$.
\end{remark}

We are going to apply the stability estimates of Theorem \ref{thm:projectors} to obtain stability estimates for eigenfunctions. For this we shall need the following lemma.
\begin{lemma}
{\bf (selection lemma)} Let $U$ and $V$ be finite dimensional subspaces of a Hilbert space $\hil$, $\dim U=\dim V=m$, and let $u_1,\ldots ,u_m$ be an orthonormal basis of $U$. Then there exists an orthonormal basis $v_1,\ldots ,v_m$ of $V$ such that
\begin{equation}
\| u_k-v_k\| \leq  5^k\|P_U-P_V\| \, , \;\; k=1,\ldots ,m \, ,
\label{basis}
\end{equation}
where $P_U$, $P_V$ are the orthogonal projectors onto $U$, $V$ respectively.
\label{burenkov}
\end{lemma}
{\em Proof.} {\em Step 1.} Clearly $\| P_U-P_V\|\le 2$. If $1\le \| P_U-P_V\| \le 2$ then estimate (\ref{basis}) obviously holds for any choice of an orthonormal basis $v_1, \dots , v_m$ of $V$ so we assume  that $\|P_U-P_V\|<1$. Let $u\in U$, $\|u\|=1$. Then
\begin{equation}
\|P_Vu\| =\| u+(P_V-P_U)u\| \geq 1-\|P_U-P_V\| >0\, .
\label{bur1}
\end{equation}
Letting $z=P_Vu / \|P_Vu\|$ we have $\|z\|=1$ and
\[
\inprod{u}{z} =\frac{\inprod{u}{P^2_Vu}}{\|P_Vu\|} =\|P_Vu\| \, ,
\]
hence
\[
\|P_U-P_V\|^2\geq \|(P_U-P_V)u\|^2 =\|u\|^2 -\|P_Vu\|^2 \geq \frac{\|u-z\|^2}{2},
\]
and therefore
\begin{equation}
\|u-z\| \leq\sqrt{2}\|P_U-P_V\|.
\label{bur2}
\end{equation}
{\em Step 2.} Assume that $\|P_U-P_V\|\le 1/6$, and
\[
z_k=\frac{P_Vu_k}{\|P_Vu_k\|} \; , \; k=1,\ldots, m.
\]
We shall prove that
\begin{equation}
|\inprod{z_k}{z_l}|\leq 3\|P_U-P_V\| \, , \ \ \;  k,l=1,\ldots, m , \; k\neq l.
\label{bur4}
\end{equation}
Indeed, we have for $k\ne l$
\begin{eqnarray*}
|\inprod{P_Vu_k}{P_Vu_l}| &=& | \inprod{P_Vu_k-u_k}{P_Vu_l} +\inprod{u_k}{u_l}+\inprod{u_k}{P_Vu_l-u_l}| \\
&=& | \inprod{ (P_V-P_U)u_k}{P_Vu_l} +\inprod{u_k}{(P_V-P_U)u_l}| \\
&\leq& 2\|P_U-P_V\|,
\end{eqnarray*}
and the claim is proved by recalling (\ref{bur1}).

{\em Step 3.} One can easily see that since $\|P_U-P_V\| <1$, the vectors $z_1,\ldots,z_m$ are linearly independent. Thus we can apply the Gram-Schmidt orthogonalization procedure, {\it i.e.},  define
\[
v_1=z_1 \; , \;\; v_k=\frac{z_k-\sum_{l=1}^{k-1}\inprod{z_k}{v_l}v_l}{\|z_k-\sum_{l=1}^{k-1}\inprod{z_k}{v_l}v_l\|} ,\ \ \  k=2,\ldots ,m.
\]
Note that for $k=2,\ldots,m$,
\[
v_k-z_k =\left(\frac{1}{\|z_k-\sum_{l=1}^{k-1}\inprod{z_k}{v_l}v_l\|}-1\right)z_k -
\frac{\sum_{l=1}^{k-1}\inprod{z_k}{v_l}v_l}{\|z_k-\sum_{l=1}^{k-1}\inprod{z_k}{v_l}v_l\|}
\]
and
\[ 1\geq
\Bigl\|z_k-\sum_{l=1}^{k-1}\inprod{z_k}{v_l}v_l\Bigr\| \geq 1-\sum_{l=1}^{k-1}|\inprod{z_k}{v_l}| \; .
\]
Hence if
\begin{equation}
\sum_{l=1}^{k-1}|\inprod{z_k}{v_l}|<1
\label{bur6}
\end{equation}
then
\begin{equation}
\|v_k-z_k\| \leq \frac{2\sum_{l=1}^{k-1}|\inprod{z_k}{v_l}|}{1-\sum_{l=1}^{k-1}|\inprod{z_k}{v_l}|}.
\label{bur7}
\end{equation}
Also for $s=k,\ldots,m$,
\begin{equation}
|\inprod{z_s}{v_k}| \leq \frac{|\inprod{z_s}{z_k}| +\sum_{l=1}^{k-1}|\inprod{z_k}{v_l}|}{1-\sum_{l=1}^{k-1}|\inprod{z_k}{v_l}|}.
\label{bur8}
\end{equation}

{\em  Step 4.}  We shall prove that for all $k=2, \dots , m$
\begin{equation}
\|v_k-z_k\| \leq 3\cdot 5^{k-1}\| P_U-P_V \| ,
\label{bur9}
\end{equation}
\begin{equation}|
\inprod{z_s}{v_k}| \leq 3\cdot 5^{k-1}\| P_U-P_V \| ,\ \ \ s=k+1 , \dots , m,
\label{bur10}
\end{equation}
provided that
\begin{equation}
\| P_U-P_V \| \le \frac{2}{3} 5^{-k+1}\, .
\label{bur11}
\end{equation}
We prove this by induction. If $k=2$ then by (\ref{bur4}) and (\ref{bur11})
$|\inprod{z_2}{ v_1}|=|\inprod{z_2}{ z_1}|\le 3 \| P_U-P_V  \|\le \frac{2}{5}  $, hence
by (\ref{bur7}),
$$
\|  v_2-z_2 \|\le \frac{6\|P_U-P_V \|}{1-3\|  P_U-P_V\| }\le 15\| P_U-P_V\|
$$
and by (\ref{bur8}) and (\ref{bur4}) for $s=3, \dots , m$ also
$$
|\inprod{z_s}{v_2}|   \le \frac{6\|P_U-P_V \|}{1-3 \|  P_U-P_V\| }\le 15\| P_U-P_V\| .
$$

Let $2\le k \le m-1$.  Assume that inequalities (\ref{bur9}) and (\ref{bur10}) under assumption (\ref{bur11}) are satisfied
for  all $2\le j \le k $. By assuming the validity of (\ref{bur11}) for $k+1$,  by (\ref{bur7}) we obtain
$$
\| v_{k+1}-z_{k+1}  \|\le \frac{6(\sum_{j=1}^k 5^{j-1} )\|  P_U-P_V\|  }{ 1-3(\sum_{j=1}^k 5^{j-1} )\|  P_U-P_V\|  } .
$$

Since $\sum_{l=1}^k5^{l-1}\le 5^k/4$, by (\ref{bur11}) with $k$ replaced by $k+1$, $3 \sum_{l=1}^k5^{l-1}\| P_U-P_V\|\le 1/2$ hence
$\| v_{k+1}-z_{k+1} \| \le 3\cdot 5^k \|  P_U-P_V\|  $.
Similarly, by (\ref{bur8}) and (\ref{bur4}) for all $s=k+2, \dots m$
$$
|(z_s, v_{k+1})|\le \frac{3 \|  P_U-P_V\| +  3(\sum_{j=1}^k 5^{j-1} )\|  P_U-P_V\|  }{ 1-3(\sum_{j=1}^k 5^{j-1} )\|  P_U-P_V\|  } \le 3\cdot 5^k \| P_U-P_V\| .
$$

{\em Step 5. } To complete the proof, we note that by (\ref{bur2}) we have $\|u_1-v_1\|\leq\sqrt{2}\|P_U-P_V\|$.
For $k\geq 2$ we have that if (\ref{bur11}) holds then
\[
\|u_k-v_k\| \leq \|u_k-z_k    \| +\|z_k -v_k\| \leq (\sqrt{2}+3\cdot 5^{k-1})\|P_U-P_V\|\, ,
\]
while if (\ref{bur11}) does not hold  then $\|P_U-P_V\|> 10/ (3\cdot 5^{k})$ and therefore
\[
\|u_k-v_k\| \leq 2 \leq 3\cdot 5^{k-1}\|P_U-P_V\|.
\]
This completes the proof of the lemma.  \hfill $\Box$

\

\begin{lemma}
\label{eigenf} Let {\rm (A)} be satisfied.
Let $\lambda$ be a non-zero eigenvalue of $H$ of multiplicity $m$ and let $n\in\N$ be such that $\lambda=\lambda_n[H]
=\ldots =\lambda_{n+m-1}[H]$. Then the following statements hold:
\begin{itemize}

\item[(i)] There exists $c_1>0$ depending only on $N,\tau ,\theta, \alpha , c^* , \lambda_{n-1}[H], \lambda , \lambda_{n+m}[H]$ such that the following is true: if
$
\delta_{\infty  }(\phi ,\tilde\phi )  \leq c_1^{-1}
$
and $\psi_n[\tilde H], \dots , \psi_{n+m-1}[\tilde H]$ are orthonormal  eigenfunctions of $\tilde H$ in $L^2(\Omega , \tilde gdx)$, then there exist orthonormal eigenfunctions
 $ \psi_n[H], \dots , \psi_{n+m-1}[H]$ of $H$ in $L^2(\Omega , gdx)$ such that
 \begin{equation}
 \label{eigenf2}
 \| \psi_{k}[H]- \psi_{k}[\tilde H]  \|_{L^2(\Omega)} \le c_1  \delta_{\infty  }(\phi ,\tilde\phi ) ,
 \end{equation}
 for all $k=n, \dots , n+m-1$.

\item[(ii)]

Let in addition {\rm (P)} be satisfied by the operators $H$, $\tilde H$ and $T^*ST$ for the same $q_0$, $\gamma $ and $C$. Let  $s=[q_0/(q_0-2)] \max\{ 2, \alpha+2\gamma \}$. Then there exists  $c_2>0$ depending only on $N,\tau ,\theta,\alpha, c^*, q_0, C, \gamma ,  |\Omega|, \lambda_{n-1}[H], \lambda ,$ $\lambda_{n+m}[H] $ such that the following is true: if
$
\delta_{s }(\phi ,\tilde\phi )  \leq c_2^{-1}
$
and $\psi_n[\tilde H]$, $\dots$ , $\psi_{n+m-1}[\tilde H]$ are orthonormal  eigenfunctions of $\tilde H$ in $L^2(\Omega , \tilde gdx)$, then there exist orthonormal eigenfunctions
$ \psi_n[H], \dots , \psi_{n+m-1}[H]$ of $H$ in $L^2(\Omega , gdx)$ such that
 \begin{equation}
 \label{eigenf4}
 \| \psi_{k}[H]-\psi_{k}[\tilde H]  \|_{L^2(\Omega)} \le c_2  \delta_{s }(\phi ,\tilde\phi ) ,
 \end{equation}
for all $k=n, \dots , n+m-1$.
\end{itemize}
\end{lemma}
{\em Proof.} We shall only prove statement $(ii)$ since the proof of statement $(i)$ is similar. We first note that $\varphi_k:=w^{-1}\psi_k[\tilde H]$, $k=n,\ldots, n+m-1$, are orthonormal eigenfunctions in $L^2(\Omega, g\, dx)$ of $w^{-1}\tilde Hw$ corresponding to the eigenvalues $\lambda _n[\tilde H], \dots , \lambda_{n+m-1}[\tilde H]$. By Theorem \ref{thm:projectors} and Lemma \ref{burenkov}
there exists $c>0$ such that if
$ \delta_s(\phi , \tilde \phi )<c^{-1}$
 then there exist eigenfunctions
$ \psi_n[H], \dots , \psi_{n+m-1}[H]$ of $H$ corresponding to the eigenvalue $\lambda$ such that
\begin{equation}
 \| \psi_{k}[H]-\varphi_{k} \|_{L^2(\Omega)} \le c  \delta_{s}(\phi ,\tilde\phi ) \, .
\label{eigenf5}
\end{equation}
In order to complete the proof it is enough to observe that
\[
\| \varphi_k -\psi_k[\tilde H]\|_{L^2(\Omega)} \leq \|1-w^{-1}\|_{L^s(\Omega)} \|\psi_k[\tilde H]\|_{L^{2s/(s-2)}(\Omega)} \leq
c \|\nabla\phi -\nabla\tilde\phi\|_{L^s(\Omega)} \, .
\]
\hfill $\Box$

\

In the following theorem we estimate the deviation of the eigenfunctions $\psi_k[\tilde L]$ of $\tilde L$ from the
eigenfunctions $\psi_k[L]$ of $L$. {\em We adopt the convention that $\psi_k[L]$ and $\psi_k[\tilde L]$ are extended by zero outside $\phi(\Omega)$ and $\tilde\phi(\Omega)$ respectively.}
\begin{theorem} {\bf (stability of eigenfunctions)} Let {\rm (A)} be satisfied.
Let $\lambda$ be a non-zero eigenvalue of $L$ of multiplicity $m$ and let $n\in\N$ be such that $\lambda=\lambda_n[L]
=\ldots =\lambda_{n+m-1}[L]$.
Then the following statements hold:
\begin{itemize}
\item[(i)]
There exists $c_1>0$ depending only on $N,\tau,\theta ,\alpha ,  c^* ,  \lambda_{n-1}[L] , \lambda ,\lambda_{n+m}[L]$ such that the following is true: if
$
\delta_{\infty  }(\phi ,\tilde\phi )\leq c_1^{-1}
$
and $\psi_n[\tilde L], \dots , \psi_{n+m-1}[\tilde L]$ are orthonormal  eigenfunctions of $\tilde L$  in $L^2(\tilde \phi (\Omega ))$, then there exist orthonormal eigenfunctions
$ \psi_n[L], \dots , \psi_{n+m-1}[L]$ of $L$ in $L^2(\phi (\Omega ) )$ such that
\begin{eqnarray}
\label{holder2}\lefteqn{
 \| \psi_{k}[L]-\psi_{k}[\tilde L]  \|_{L^2(\phi (\Omega )\cup \tilde \phi (\Omega ))} \le  c \big(
 \delta_{\infty }(\phi ,\tilde\phi ) +}\nonumber \\
& &
+  \| \psi_k[L]\circ \phi -\psi_k[L]\circ \tilde \phi  \|_{L^{2 }(\Omega ) }+
  \| \psi_k[\tilde L]\circ \phi -\psi_k[\tilde L]\circ \tilde \phi  \|_{L^{2 }(\Omega ) }  \big)\! ,
 \end{eqnarray}
 for all $k=n, \dots , n+m-1$.
\item[(ii)]
Let in addition {\rm (P)} be satisfied by the operators $L$, $\tilde L$ and $\hat L$ for the same $q_0$, $\gamma$ and $C$.
Let $s=[q_0/(q_0-2)] \max\{ 2, \alpha+2\gamma \}$. Then
there exists  $c_2>0$ depending only on $N,\tau ,\theta,\alpha, c^*, q_0, C, \gamma ,  |\Omega|,\lambda_{n-1}[L] , \lambda $, $ \lambda_{n+m}[L]  $, such that the following is true: if
$
\delta_{s  }(\phi ,\tilde\phi )\leq c_1^{-1}
$
and $\psi_n[\tilde L], \dots $, $\psi_{n+m-1}[\tilde L]$ are orthonormal  eigenfunctions of $\tilde L$  in $L^2(\tilde \phi (\Omega ))$, then there exist orthonormal eigenfunctions
$ \psi_n[L], \dots , \psi_{n+m-1}[L]$ of $L$ in $L^2(\phi (\Omega ) )$ such that
\begin{eqnarray}
\label{holder2bis}\lefteqn{
 \| \psi_{k}[L]-\psi_{k}[\tilde L]  \|_{L^2(\phi (\Omega )\cup \tilde \phi (\Omega ))} \le  c \big(
 \delta_{s }(\phi ,\tilde\phi ) +}\nonumber \\
& &
+  \| \psi_k[L]\circ \phi -\psi_k[L]\circ \tilde \phi  \|_{L^{2 }(\Omega ) }+
  \| \psi_k[\tilde L]\circ \phi -\psi_k[\tilde L]\circ \tilde \phi  \|_{L^{2 }(\Omega ) }  \big)\! ,
 \end{eqnarray}
 for all $k=n, \dots , n+m-1$.
\end{itemize}
\label{holder}
\end{theorem}
\begin{remark} We note that if in addition the semigroup $e^{-Lt}$ is ultracontractive then
the eigenfunctions are bounded hence
$$\| \psi_k[L]\circ \phi -\psi_k[L]\circ \tilde \phi  \|_{L^{2 }(\Omega )}+\| \psi_k[\tilde L]\circ \phi -
\psi_k[\tilde L]\circ \tilde \phi  \|_{L^{2 }(\Omega )}\leq c(\lambda) |\cD|^{1/2},
$$
where $\cD=\{x\in \Omega :\, \phi (x)\ne \tilde\phi (x)  \}$.
\end{remark}
{\em Proof of Theorem \ref{holder}.} We set $$\psi_k[\tilde H] =\psi_k[\tilde L] \circ \tilde \phi , $$ for all $k=n, \dots , n+m-1$, so that $\psi_n[\tilde H], \dots ,$ $\tilde \psi_{n+m-1}[\tilde H]$ are orthonormal eigenfunctions in $L^2(\Omega , \tilde g dx)$ of the operator $\tilde H$ corresponding to the eigenvalues $\lambda_n[\tilde H], \dots ,$ $ \lambda_{n+m-1}[\tilde H]$. By Lemma~\ref{eigenf} $(i)$ it follows that there exists $c_1>0$ such that if $\delta_{\infty }(\phi , \tilde \phi)<c_1^{-1}$ then there exist  orthonormal  eigenfunctions $ \psi_n[H], \dots ,$ $ \psi_{n+m-1}[H]$ in $L^2(\Omega ,  g dx)$ of $H$ corresponding to the eigenvalue $\lambda$ such that inequality (\ref{eigenf2}) is satisfied. We now set
$$
\psi_k[L]= \psi_k[H]\circ \ \phi^{(-1)},
$$
for all $k=n, \dots , n+m-1$, so that $ \psi_n[L], \dots , \psi_{n+m-1}[L]$ are orthonormal eigenfunctions  in $L^2(\phi (\Omega ) )$ of $L$ corresponding to the eigenvalue $\lambda$.
Then by changing variables in integrals we obtain
\begin{eqnarray}\lefteqn{
 \|\psi_k[L]-\psi_k[\tilde L] \|_{L^2(\tilde\phi(\Omega))} \le
  \|\psi_k[L]\circ\tilde \phi-\psi_k[\tilde L]\circ\tilde\phi \|_{L^2(\Omega)} } \nonumber \\ & & \qquad\qquad
 \le  c \Big( \|\psi_k[L]\circ\tilde\phi-\psi_k[L]\circ\phi \|_{L^2(\Omega)} +\|\psi_k[L]\circ \phi-\psi_k[\tilde L]\circ\tilde\phi \|_{L^2(\Omega)}\Big) \nonumber \\
 & &  \qquad\qquad  = c \Big( \|\psi_k[L]\circ\tilde\phi-\psi_k[L]\circ\phi \|_{L^2(\Omega)} +\|\psi_k[H]-\psi_k[\tilde H]\|_{L^2(\Omega)}\Big) . \nonumber
\end{eqnarray}
In the same way $$\|\psi_k[L]-\psi_k[\tilde L] \|_{L^2(\phi(\Omega))} \le c(
\|\psi_k[\tilde L]\circ\tilde\phi-\psi_k[\tilde L]\circ\phi \|_{L^2(\Omega)} +\|\psi_k[H]-\psi_k[\tilde H]\|_{L^2(\Omega)}).$$
Hence estimates (\ref{holder2}), (\ref{holder2bis}) follow by (\ref{eigenf2}), (\ref{eigenf4}) respectively. \hfill $\Box$

\section{On regularity of eigenfunctions}

\label{reg}

In this section we  obtain sufficient conditions for the validity of conditions (P1) and (P2).
We begin by recalling the following known  result based on the notion of ultracontractivity which guarantees the validity of property (P1) under rather general assumptions, namely under the assumption that a Sobolev-type  Embedding Theorem  holds for the space ${\mathcal{V}}$.

\begin{lemma}
\label{ultra}
Let $\Omega $ be a domain in $\R^N$ of finite measure and ${\mathcal{V}}$ a closed subspace of $W^{1,2}(\Omega )  $ containing $W^{1,2}_0(\Omega )  $.
Assume that there exist $p>2, D>0$ such that
\begin{equation}
\label{sobolev}
\| u\| _{L^p(\Omega )}\le D \| u\| _{W^{1,2}(\Omega )} ,
\end{equation}
for all $u\in \cV$. Then the following statements hold:
\begin{itemize}
\item[(i)] Condition (\ref{cstar}) is satisfied for any $\alpha>\frac{p}{p-2}$.
\item[(ii)] The  eigenfunctions of the operators $H$, $\tilde H$ and $T^*ST$ satisfy  {\rm (P1)} with $q_0=\infty$, $\gamma =\frac{p}{2(p-2)}$,  where $C$ depends only on $p, D, \tau,\theta , c^* $.
\end{itemize}
\end{lemma}
{\em Proof.} For the proof of statement $(i)$ we refer to \cite[Thm.~7]{buda} where the case ${\mathcal{V}}=W^{1,2}(\Omega )$ is considered. The proof works word by word also in the slightly more general case considered here. The proof of statement $(ii)$ is as in \cite[Thm.~7]{buda} where it is proved that for the Neumann Laplacian property (P1) is satisfied if (\ref{sobolev}) holds: this proof can be easily adapted to the operators $H$, $\tilde H$ and $T^*ST$. \hfill $\Box$\\

We now give conditions for the validity of property (P2). We consider first the case when an  a priori estimate holds for the operators $L$, $\tilde L$, which is typically the case of sufficiently smooth open sets and coefficients. Then we consider a more general situation based on an approach which goes back to Meyers~\cite{mey}.

\

{\bf The regular case}

\

Recall that an open set in $\R^N$ satisfies the interior cone
condition with the parameters $R>0$ and $h>0$ if for all
$x\in\Omega$ there exists a cone $K_x\subset \Omega $ with the point $x$ as
vertex congruent to the cone
$$
K(R,h)=\biggl\{x\in\R^N:\
0<\biggl(\sum_{i=1}^{N-1}x_i^2\biggr)^{1/2} <\frac{Rx_N}{h} <R
\biggr\}.
$$

In this paper the cone condition is used in order to guarantee the validity of the standard Sobolev embedding.

The next theorem is a simplified version of \cite[Theorem~5.1]{bulaneu}.
\begin{theorem}
\label{apriori0} Let $R>0$,
$h>0$. Let $U$ be an open set in $\R^N$
satisfying the interior cone condition with the parameters $R$ and $h$,
and let $E$ be an operator in  $L^{2}(U)$ satisfying the following
a priori estimate:
\begin{enumerate}
\item[] there exists $B>0$ such that if  $2\le p <N+2$ \hfill and if
\hfill $u\in {\rm Dom}(E)$ \hfill
 and\hfill\hfill\hfill\hfill\hfill\\
$Eu\in L^p(U)$, then $u\in W^{2,p}(U)$ and
\begin{equation}
\label{apriori} \|u \|_{  W^{2,p}(U) }  \le B \left( \|Eu
\|_{L^p(U)} +\|u \|_{  L^2(U) }  \right).
\end{equation}
\end{enumerate}

Assume that $E\psi =\lambda \psi $ for some $\psi \in {\rm Dom}(E)$ and $\lambda\in
\C$.
Then
 there exists $c>0$, depending only on $
  R, h, N$ and $B$, such that for $\mu =0,1$,
\begin{equation}
\label{apriori1}  \|\psi \|_{  W^{\mu ,\infty}(U) } \le
c(1+|\lambda |)^{\frac{N}{4}+
\frac{\mu }{2}}
\|\psi \|_{ L^{2}(U) }.
\end{equation}
\end{theorem}

\begin{theorem}
\label{troy}
Let {\rm (A)} be satisfied and let
$\phi(\Omega )$ and $\tilde \phi(\Omega )$ be open sets satisfying the interior cone condition with the same parameters $R,h$. If
the operators $L$, $\tilde L$ satisfy the a priori estimate (\ref{apriori}) with the same $B$, then
the operators $H$, $\tilde H$, $T^*ST$ satisfy property {\rm (P)} with $q_0=\infty $, $\gamma =N/4$ and $C$ depending only on $\tau , R, h, c^*, \theta $ and $B$.
\end{theorem}
{\em Proof. } Recall that $H, \tilde H$ and $T^*ST$ are the operators obtained by pulling-back to $\Omega$ the operators $L, \tilde L$ and $\hat L$ respectively. Clearly $\hat L$ also satisfies the a priori estimate (\ref{apriori}). Thus, by Theorem~\ref{apriori0} the eigenfunctions
of the operators $L, \tilde L, \hat L$ satisfy condition (\ref{apriori1}) hence, by pulling such eigenfunctions back to $\Omega $ it follows that
the eigenfunctions of $H$, $\tilde H$, $T^*ST$ satisfy (P1) and (P2) with  $q_0=\infty $, $\gamma =N/4$ and $C$ as in the statement.
\hfill $\Box$

\

{\bf The general case}

\

Here we shall assume that ${\mathcal{V}}= {\rm cl}_{W^{1,2}(\Omega)}\cV_0$ where $\cV_0$
is a space of functions defined on $\Omega $ such that $C_c^{\infty}(\Omega)\subset \cV_0\subset  W^{1,\infty}(\Omega)$. Moreover, for all $1 < q  < \infty $ we set
$$
V_q= {\rm cl}_{W^{1,q}(\Omega)}\cV_0.
$$

Let  $-\Delta_{q}: V_q\to (V_{q'})'$ be the operator defined by
$$
(-\Delta_{q}u, \psi)=\int _{\Omega }\nabla u\cdot \nabla \psi dx,
$$
for all $u\in V_q$, $\psi\in V_{q'}$.

The following theorem is a variant of a result of Gr\"oger \cite{groger}; see also \cite{bar3}.

\begin{theorem}\label{meyer}
Let {\rm  (A)} be satisfied. Assume that
there exists $q_1>2$ such that
the operator $I-\Delta_{q} : V_q\to (V_{q'})'$ has a bounded inverse for all $2\le q\le q_1$.
Then there exist $q_0>2$ and $c>0$, depending only on $\cV_0$, $\tau$ and $\theta$ such that if $u$ is an eigenfunction of one  of the operators $H$, $\tilde H$, $T^*ST$ and $\lambda $ is the corresponding eigenvalue then
\begin{equation}
\label{mey1}
\| \nabla u\|_q\le c(1+\lambda )\| u \| _q,
\end{equation}
for all $2\le q\le q_0$.

Moreover, if $\Omega$ is such that the interior cone condition holds then there exists $c>0$ depending only on
$\cV$, $\tau$ and $\theta$ such that
\begin{equation}
\label{mey2}
\| \nabla u\|_q\le c(1+\lambda )\| u \| _{\frac{Nq}{N+q}},
\end{equation}
for all  $2< q\le q_0$.
\end{theorem}
{\em Proof.} We prove the statement for the operator $T^*ST$, the other cases being similar. We divide the proof into three steps.

{\it Step 1.}  We define
$$
Q(u,\psi )=\int_{\Omega } u\psi gdx +\int_{\Omega }\tilde a \nabla u\cdot \nabla \psi \tilde gdx,
$$
$$
Q_0(u,\psi )=\int_{\Omega } u \psi dx +\int_{\Omega } \nabla u\cdot \nabla \psi dx,
$$
for all $u\in V_q$, $\psi \in V_{q'}$. Since\footnote{
Here we use
$\|f\|_{W^{1,p}(\Omega )}^p=\|f\|_{L^p(\Omega )}^p +\|   \, |\nabla f| \, \|_{L^p(\Omega )}^p$ as the norm in $W^{1,p}(\Omega ) $.}
$$
|Q_0(u,\psi ) -\beta Q(u,\psi ) |\le \max\{ \| 1-\beta g \|_{L^{\infty}(\Omega)},
\|I- \beta\tilde a\tilde g\|_{L^{\infty}(\Omega)}\} \| u\|_{W^{1,q}(\Omega ) }\| \psi \|_{W^{1,q'}(\Omega ) },
$$
there exist $\beta>0$ and $0<c<1$ depending only on $N ,\tau$ and $\theta$ such that
\begin{equation}
\label{quadest}
|Q_0(u,\psi ) -\beta Q(u,\psi )  |\le c \|u\|_{W^{1,q}(\Omega ) }\| \psi \|_{W^{1,q'}(\Omega ) },
\end{equation}
for all $u\in W^{1,q}(\Omega)$ and $\psi \in W^{1,q'}(\Omega )$.

{\it Step 2.}  Using the fact that $\| (I-\Delta_2)^{-1}  \|=1$, that $q\mapsto \| (I-\Delta_q)^{-1}\|$
is continuous and by observing that $2 /(c+1)>1$, it follows that there exists $q_0>2$ such that
\begin{equation}
\label{quadest1}
\| (I-\Delta_q)^{-1}\| <\frac{2}{c+1},
\end{equation}
for all $2\le q \le q_0$. By (\ref{quadest}) it then follows that for all $2\leq q\leq q_0$,
\begin{eqnarray}
\inf_{\| u\|_{W^{1,q}(\Omega )}=1}\sup_{\|\psi \|_{W^{1,q'}(\Omega )}=1} Q(u,\psi )&\geq&
\frac{1}{\beta}\inf_{\|\psi\|_{W^{1,q'}(\Omega )}=1}\sup_{\|u\|_{W^{1,q}(\Omega )}=1}Q_0(u,\psi )-\frac{c}{\beta} \nonumber \\
&=&\frac{1}{\beta}\| (I-\Delta_q)^{-1}  \|^{-1}-\frac{c}{\beta}\nonumber \\
&>& \frac{1-c}{2\beta}>0\, .
\label{quadest1.5}
\end{eqnarray}
{\it Step 3.} By (\ref{quadest1.5}) it follows that the operator $I+   (T^*ST)_q$ of $V_q$ to $V_{q'} '$ defined by
\begin{equation}
\label{quadest2}
(I+  (T^*ST)_qu,\psi )=Q (u, \psi )
\end{equation}
has a bounded inverse such that
\begin{equation}
\label{quadest2bis}
\|( I+  (T^*ST)_q)^{-1}\| = \biggl(\inf_{\| u \|_{W^{1,q}(\Omega )} =1}\sup_{\| \psi \|_{W^{1,q'}(\Omega )}=1} Q (u,\psi )\biggr)^{-1}<\frac{2\beta} {1-c}.
\end{equation}
Then (\ref{mey1}) follows by (\ref{quadest2}), (\ref{quadest2bis}) and by observing that
\begin{equation}
\label{weakreg}
Q (u,\psi )=(1 + \lambda )\int _{\Omega }u\psi g\,  dx,
\end{equation}
for all $\psi \in V_{q'}$.

Now, if $\Omega$ satisfies the interior cone condition, then the standard Sobolev embedding holds. Thus, if $q>2$ then $q'<2\le N$, hence $V_{q'}$ is continuously embedded into $L^{\frac{Nq'}{N-q'}}(\Omega )$.
By (\ref{weakreg})
we have
\begin{eqnarray}
\| u\|_{W^{1,q}(\Omega )} & \le &  (1 +\lambda ) \| ( I+ (T^*ST)_q)^{(-1)}  \|\sup_{\|\psi\|_{W^{1,q'}(\Omega )}=1}\Big | \int_{\Omega }u\psi g\, dx\Big |\nonumber \\
& \le & \frac{2\beta }{1-c}(1 +\lambda )\| g\|_{L^{\infty}(\Omega ) }\| u \|_{L^{\frac{Nq}{N+q}}(\Omega )} \sup_{\|\psi\|_{ W^{1,q'}(\Omega )  }=1}
\| \psi \|_{L^{\frac{Nq'}{N-q'}}(\Omega )} \, ,
\end{eqnarray}
and the last supremum is finite due to the Sobolev embedding.
\hfill $\Box$

\

\begin{remark}
\label{ultrameyer}
If $\Omega$ satisfies the interior cone condition then inequality (\ref{sobolev}) is satisfied with $p=2N/(N-2)$ if $N\geq 3$ and with any $p>2$ if $N=2$. Then  by Lemma~\ref{ultra} it follows that condition (\ref{cstar})
holds for any $\alpha >N/2$ and
the operators $H, \tilde H, T^*ST, L, \tilde L, \hat L$ satisfy property {\rm(P1)} with $q_0 =\infty $,  $\gamma = N/4$ if $N\geq 3$ and any $\gamma >1/2$ if $N=2$.
In fact, if $N=2$ property {\rm(P1)} is also satisfied for $\gamma =1/2$; this follows by \cite[Thm.~2.4.4]{daheat} and \cite[Lemma~10]{buda}.
Thus by the second part of Theorem~\ref{meyer} it follows that both properties (P1) and (P2) are satisfied for some $q_0>2$ and $\gamma = N(q_0-2)/(4q_0)$ for any $N\geq 2$.

If $\Omega$ is of class $C^{0,\nu }$ ({\it i.e.,} $\Omega$ is locally a subgraph of $C^{0,\nu}$ functions)
with $0<\nu <1$, then inequality (\ref{sobolev}) is satisfied with $p=2(N+\nu-1)/(N-\nu-1)$, for any $N\geq 2$ (see also \cite{buda}).
Thus Lemma \ref{ultra} implies that condition (\ref{cstar})
holds for any $\alpha > (N+\nu -1)/(2\nu)$
and that the operators $H, \tilde H, T^*ST, L, \tilde L, \hat L$ satisfy property {\rm(P1)} with $q_0 =\infty $ and  $\gamma =(N+\nu-1)/(4\nu)$.
\end{remark}

\section{Estimates via Lebesgue measure}\label{applsec}

In this section we consider two examples to which we apply the results of the previous sections in order to obtain stability estimates via the Lebesge measure.

Let $A_{ij}\in L^{\infty}(\R^N)$ be real-valued functions satisfying $A_{ij}=A_{ji}$ for all $i,j=1, \dots , N$ and condition (\ref{ellip}). Let $\Omega  $ be a bounded domain in $\R^N$ of class $C^{0,1}$,
and let $\Gamma$ be an open subset of $\partial \Omega$ with a Lipschitz boundary in    $\partial \Omega $ (see Definition~\ref{lipgammabd} below).  We consider the eigenvalue problem with mixed Dirichlet-Neumann boundary conditions
\begin{equation}
\label{mixed}
\left\{
\begin{array}{ll}
-\sum_{i,j=1}^N\frac{\partial}{\partial x_i }\big( A_{ij}(x)\frac{\partial u}{\partial x_j} \big)=\lambda u, &\ {\rm in}\ \Omega ,
\vspace{0.2cm}\\
u=0, &\ {\rm on}\ \Gamma ,\vspace{0.2cm}\\
\sum_{i,j=1}^N A_{ij}\frac{\partial u}{\partial x_j}\nu_{i}=0, &\ {\rm on}\ \partial\Omega  \setminus \Gamma ,
\end{array}\right.
\end{equation}
where $\nu$ denotes the exterior unit normal to $\partial \Omega$. Observe that our analysis comprehends the `simpler' cases   $\Gamma =\partial \Omega   $ (Dirichlet boundary conditions) or $\Gamma =\emptyset$ (Neumann boundary conditions), as well as all other cases where
$\Gamma $ is a connected component of $\partial \Omega$ (the boundary of $\Gamma$ in $\partial \Omega $ is empty).  See \cite{groger} for details.

We denote by $\lambda_{n}[\Omega ,\Gamma ]$ the sequence of the eigenvalues of problem (\ref{mixed}) and by   $\psi_n [\Omega ,\Gamma]  $ a corresponding orthonormal system of  eigenfunctions in $L^2(\Omega )$. In this section we  compare the eigenvalues and the eigenfunctions corresponding to open sets $\Omega$ and $\tilde\Omega$ and
the associate portions of the boundaries $\Gamma\subset\partial \Omega $ and $\tilde \Gamma\subset \partial \tilde\Omega$.
To do so we shall think of $\Omega $ as a fixed reference domain and we shall apply the results of the previous sections
to transformations $\phi$ and $\tilde\phi $ defined on $\Omega$, where $\phi =Id $ and $\tilde \phi$ is a suitably constructed bi-Lipschitz homeomorphism such that
$\tilde \Omega =\tilde\phi (\Omega )$ and $\tilde \Gamma =\tilde \phi (\Gamma )$.

Before doing so, we recall the weak formulation of problem (\ref{mixed}) on $\Omega $. Given $\Gamma \subset \partial \Omega $ we consider the space
$
 W^{1,2}_{\Gamma }(\Omega )
$
 obtained by taking the closure of $ C^{\infty }_{\Gamma }(\overline{ \Omega }) $ in $W^{1,2}( \Omega )$, where $ C^{\infty }_{\Gamma }(\overline{ \Omega }) $ denotes the space of the functions in $ C^{\infty }(\overline{\Omega   }) $
 which vanish in a neighborhood of $\Gamma$. Then the eigenvalues and eigenfunctions
of problem (\ref{mixed}) on $\Omega $ are the eigenvalues and the eigenfunctions of the operator $L$ associated with the sesquilinear form   $Q_{L}$  defined on ${\mathcal{W}}:= W^{1,2}_{\Gamma }( \Omega )$ as in (\ref{quadratic}).

\begin{definition}
\label{lipgammabd}
Let $\Omega$ be a bounded  open set in $\R^N$ of class $C^{0,1}$ and let $\Gamma$ be an open subset of $\partial \Omega$. We say that $\Gamma$ has a Lipschitz continuous boundary $\partial \Gamma $ in $\partial\Omega $ if for all $x\in \partial \Gamma$ there exists an open neighborhood $U$ of $x$ in $\R^N$ and $ \phi \in \Phi (U)$
such that $$  \phi ( U\cap (\Omega \cup \Gamma))=
 \{x\in \R^N\!:\, |x|<1,\, x_N<0    \}\cup \{x\in \R^N\!:\, |x|<1,\, x_N\le 0,\,  x_1>0    \} \, .
$$
\end{definition}

\subsection{Local perturbations}
\label{localsec}

In this section we consider open sets belonging to the following class.

\begin{definition}
Let $V$ be a bounded open cylinder, i.e., there exists a rotation  $R$  such that $R(V)=W\times ]a,b[$, where $W$ is a bounded convex open set in $\R^{N-1}$. Let $M,\rho>0$. We say that a bounded open set $\Omega\subset\R^N$ belongs to $\cC^{m,1}_M(V,R,\rho)$ if
$\Omega$ is of class $C^{m,1}$ (i.e., $\Omega$ is locally a subgraph of $C^{m,1}$ functions) and there exists a function $g\in C^{m,1}(\overline{W})$ such that $a+\rho\leq g\le b$,   $| g |_{m,1}:=\sum_{0<|\alpha |\le m+1}\| D^{\alpha }g \|_{L^{\infty }(W )}\le M$,
 and
\begin{equation}
R(\Omega\cap V)=\{ (\bar{x},x_N) \; : \;  \bar{x}\in W \, ,\,  a<x_N<g(\bar{x}) \}.
\label{C11}
\end{equation}
\label{def:C11}
\end{definition}

Let $\Omega ,  \tilde \Omega \in \cC^{0,1}_M(V,R,\rho)$ be such that $\Omega \cap (V_{\rho})^{c}=\tilde \Omega \cap (V_{\rho})^{c}$.   We shall assume that the corresponding sets $\Gamma \subset \partial \Omega  $, $\tilde \Gamma \subset \partial\tilde  \Omega $, where Dirichlet boundary conditions are imposed, are such that
\begin{equation}
\label{bdmixed}
\Gamma \cap V^{c}=\tilde \Gamma \cap V^{c},\ \ {\rm and}\ \ P_{R^{(-1)}W}(\Gamma \cap V)=P_{R^{(-1)}W}(\tilde \Gamma \cap V),
\end{equation}
where $P_{R^{(-1)}W}$ denotes the orthogonal projector onto $R^{(-1)}W$.
Given $\Gamma$, condition (\ref{bdmixed}) uniquely determines $\tilde\Gamma$.

\begin{theorem}
Let $\Omega \in \cC^{0,1}_M(V,R,\rho)$ and let $\Gamma $ be an open subset of
$\partial \Omega $ with Lipschitz continuous boundary in $\partial\Omega  $. Then there exists $2<q_0\leq\infty$ such that
for any $r>\max \{ 2,N(q_0-1)/q_0\} $
the following statements hold:
\begin{itemize}
\item[(i)]
There exists $c_1>0$ such that
\begin{equation}
\left(\sum_{n=1}^{\infty}\left| \frac{1}{
{\lambda}_n[\tilde\Omega , \tilde \Gamma ] +1} - \frac{1}{{\lambda}_n[\Omega , \Gamma ] +1}\right|^r\right)^{1/r}
\leq c_1 |\Omega \vartriangle \tilde \Omega  |^{\frac{q_0-2}{rq_0}} ,
\end{equation}
for all $\tilde \Omega \in  \cC^{0,1}_M(V,R,\rho) $ such that $\tilde\Omega \cap (V_{\rho})^{c}=\Omega \cap (V_{\rho})^{c}$,
$|\Omega \vartriangle \tilde \Omega |\le c_1^{-1}
$,
where  $\tilde \Gamma \subset \partial \tilde \Omega $ is determined by condition (\ref{bdmixed}).
\item[(ii)] Let $\lambda [\Omega ,\Gamma ]$ be an eigenvalue of multiplicity $m$ and let $n\in \N$ be such that
$\lambda [\Omega ,\Gamma ] =\lambda_n[\Omega ,\Gamma ]=\dots =\lambda_{n+m-1}[\Omega , \Gamma ]$.
There exists $c_2>0$ such that the following is true: if
$\tilde \Omega \in  \cC^{0,1}_M(V,R,\rho) $, $\Omega \cap (V_{\rho})^{c}=\tilde \Omega \cap (V_{\rho})^{c}$,
$
|\Omega \vartriangle \tilde \Omega |\le c_2^{-1},
$
and  $\tilde \Gamma\subset \partial \tilde \Omega $ is determined by  (\ref{bdmixed}) then, given orthonormal eigenfunctions
$\psi_n[\tilde \Omega ,\tilde \Gamma ], \dots , $ $ \psi_{n+m-1}[\tilde \Omega ,\tilde \Gamma ]$, there exist corresponding  orthonormal eigenfunctions $ \psi_n[\Omega ,\Gamma ],$ $ \dots ,$   $  \psi_{n+m-1}[\Omega ,\Gamma ] $
such that
$$
\| \psi_n[\Omega ,\Gamma ]-\psi_n[\tilde \Omega,\tilde \Gamma ] \|_{L^2(\Omega \cup\tilde \Omega)}\le c_2|\Omega \vartriangle\tilde \Omega |^{\frac{q_0-2}{rq_0}}.
$$
\end{itemize}
Moreover, if in addition $A_{ij}\in C^{0,1}(\R^N)$, $\Omega, \tilde \Omega \in  \cC^{1,1}_M(V,R,\rho)$ and $\Gamma$ is a connected component of
$\partial \Omega $ then statements (i) and (ii) hold with $q_0=\infty $.

\label{thm:C11}
\end{theorem}

For the proof we  need the following variant of \cite[Lemma 4.1]{bula}.
\begin{lemma}
\label{graf} Let $W$ be a bounded convex open set in
$\R^{N-1}$ and $M>0$. Let $0<\rho<b-a$ and $g_1,g_2$ be Lipschitz
continuous functions from $\overline{W}$ to $\R$ such that
\begin{equation}
\label{graf1}
a+\rho < g_1(\bar x) , g_2(\bar x) <b ,
\end{equation}
for all $\bar x\in \overline{W}$, and such that ${\rm Lip}g_1$, ${\rm Lip}g_2\le M $  . Let $\delta =  \frac{\rho }{2(b-a)}    $ and
$g_3= \min \{g_1 ,g_2\}-\delta |g_1-g_2  |$. Let
\begin{equation}
\label{graf2}
{\mathcal{O}}_k:= \left\{(\bar x,x_N):\ \bar x\in W,\ a<x_N<g_k(\bar x) \right\}
\end{equation}
for $k=1,2,3$.
Let $\Phi$ be the map from  $\overline{{\mathcal{O}}}_1$ into $\overline{{\mathcal{O}}}_2$
defined as follows:
\begin{itemize}
\item[] if  $g_2(\bar x)\le g_1(\bar x)  $ then
\begin{equation}
\label{graf3bis}
\Phi (\bar x,x_N)\equiv \left\{
\begin{array}{ll}
(\bar x,x_N) & {\rm if}\ (\bar x,x_N)\in  {\overline{{\mathcal{O}}}}_{3 }\\
\left(\bar x, g_2(\bar x)+\frac{\delta}{\delta +1}(x_N-g_1(\bar x)) \right)
 & {\rm if}\ (\bar x,x_N)\in {\overline{{\mathcal{O}}}}_1\setminus   {\overline{{\mathcal{O}}}}_{3},
\end{array}
\right.
\end{equation}
\item[] while if $g_2(\bar x)> g_1(\bar x)  $ then
\begin{equation}
\label{graf3}
\Phi (\bar x,x_N)\equiv \left\{
\begin{array}{ll}
(\bar x,x_N) & {\rm if}\ (\bar x,x_N)\in  {\overline{{\mathcal{O}}}}_{3 }\\
\left(\bar x, g_2(\bar x)+\frac{\delta +1}{\delta }(x_N-g_1(\bar x)) \right)
 & {\rm if}\ (\bar x,x_N)\in {\overline{{\mathcal{O}}}}_1\setminus   {\overline{{\mathcal{O}}}}_{3}.
\end{array}
\right.
\end{equation}
\end{itemize}

Then $\emptyset \ne {\mathcal{O}}_{3}\subset {\mathcal{O}}_{1}\cap {\mathcal{O}}_{2}$,
\begin{equation}
\label{graf4}
| \{ x\in {\mathcal{O}}_1:\ \Phi (x)\ne x \} | =  |{\mathcal{O}}_{1}\setminus {\mathcal{O}}_{3}|\le 2|{\mathcal{O}}_{1}\vartriangle {\mathcal{O}}_{2}|,
\end{equation}
and
$\Phi $ is a bi-Lipschitz  homeomorphism of
$\overline{{\mathcal{O}}}_1$ onto $\overline{{\mathcal{O}}}_2$.  Moreover  $\Phi \in \Phi_{\tau}(\Omega )$ where $\tau$ depends only on $N, M, \delta $.
\end{lemma}

{\em Proof.} The proof is as in \cite[Lemma~4.1]{bula} where the case $g_2\le g_1$ was considered: here we simply replace $g_1-g_2$ by $|g_1-g_2|$. \hfill $\Box$\\

{\em Proof of Theorem \ref{thm:C11}.} We shall apply Theorems \ref{thm:series} and \ref{holder} with $\phi =Id$ and $\tilde\phi$ given by
\begin{equation}
\tilde\phi(x)=\left\{
\begin{array}{ll}
x, & x\in \overline{\Omega }\setminus V , \\
R^{(-1)}\circ \Phi \circ R (x) , &  x\in \overline{\Omega }\cap V.
\end{array}\right.
\end{equation}
Here $\Phi$ is defined as in Lemma \ref{graf} for $g_1=g$ and $g_2=\tilde g$, where $g, \tilde g$ are the functions describing
the boundary in $V$  of $\Omega$, $\tilde \Omega $ respectively, as in Definition~\ref{def:C11}. Then clearly $\phi,\tilde\phi\in\Phi_{\tau}(\Omega)$, where $\tau$ depends only on $N, V, M, \rho$.
Clearly $\phi(\Omega)=\Omega$ and $\tilde\phi(\Omega)=\tilde\Omega$.  Moreover,  $\tilde \phi (\Gamma )=\tilde \Gamma $, hence $$ C_{\tilde \phi }[W_{\tilde \Gamma}^{1,2}(\tilde \Omega )  ]= C_{ \phi }[W_{ \Gamma}^{1,2}(\Omega )  ].$$
Moreover, condition (\ref{cstar}) is satisfied for any $\alpha >N/2$, see Remark~\ref{safarov}.  Hence assumption (A) is satisfied.
Observe that by (\ref{graf4}) and by the boundedness of the coefficients $A_{ij}$,
\begin{equation}
\label{meas}
\delta_p(\phi,\tilde\phi)^p\leq c\int_{\{x\in\Omega : \phi(x)\neq\tilde\phi(x)\}} \big(|\nabla\phi-\nabla\tilde\phi|^p + |A\circ\phi-A\circ\tilde\phi|^p\big)dx\leq c|\Omega\vartriangle\tilde\Omega|.
\end{equation}
By \cite[Theorem 3]{groger} the assumption of Theorem \ref{meyer} is satisfied for the space
$\cV_0=C_{\Gamma}^{\infty}(\bar{\Omega})$ for some $2<q_1<\infty$. Thus by Remark~\ref{ultrameyer}  the operators $L$, $\tilde L$ and $\hat L$ satisfy properties (P1) and (P2) for some $2<q_0<\infty$ and $\gamma=N(q_0-2)/(4q_0)$.
Thus statement $(i)$  follows by  Theorem~\ref{thm:series}~$(ii)$ with $p=q_0/(q_0-2)$. Moreover, Theorem~ \ref{holder}~$(ii)$ provides the existence
of orthonormal eigenfunctions $\psi_k[\Omega , \Gamma]$ satisfying estimate  (\ref{holder2bis}) with  $s=[q_0/(q_0-2)]
\max\{ 2, N(q_0-1)/q_0\}$. By Lemma~\ref{ultra} the
functions $\psi_k[\Omega , \Gamma ]$, $\psi_k[\tilde \Omega ,\tilde  \Gamma ]$ are bounded, hence by (\ref{graf4})
\begin{equation}
\label{uni}
\| \psi_k[\Omega , \Gamma ]\circ \phi -\psi_k[\Omega , \Gamma]\circ \tilde \phi  \|_{L^{2 }(\Omega ) }^2 ,
\| \psi_k [\tilde \Omega ,\tilde  \Gamma] \circ \phi -\psi_k[\tilde \Omega ,\tilde  \Gamma]\circ \tilde \phi  \|_{L^{2 }(\Omega ) }^2\le c |\Omega \vartriangle \tilde\Omega | .
\end{equation}
Thus statement $(ii)$ follows by estimates (\ref{holder2bis}) and (\ref{uni}).

Finally, if $A_{ij}\in C^{0,1}(\R^N)$, $\Omega, \tilde \Omega \in  \cC^{1,1}_M(V,R,\rho)$ and $\Gamma$ is a connected component of
$\partial \Omega $, by Troianiello~\cite[Thm.~3.17~(ii)]{tro} the operators $L$ and $\tilde{L}$ satisfy the a priori estimate (\ref{apriori}) on $\Omega$ and $\tilde\Omega$ respectively. Thus by Theorem \ref{troy} the operators $L$, $\tilde L$ and $\hat L$ satisfy properties (P1) and (P2) with $q_0=\infty$ and $\gamma=N/4$, and the result follows as above.  $\hfill \Box$

\subsection{Global normal perturbations}
\label{globsec}

Let $\Omega$ be a bounded domain with $C^{2}$ boundary. By the Tubular Neighborhood Theorem  there exists $t>0$ such that for each $x\in (\partial\Omega )^t :=\{x\in \R^N \, : \;  \dist(x,\partial\Omega)<t\}$ there exists a unique couple $(\bar{x},s)\in\po\times ]-t,t[ $ such that $x=\bar{x}+s\nu(\bar{x})$; moreover, $\bar{x}$ is the (unique) nearest to $x$ point of the boundary and $s=\dist(x,\partial\Omega)$. One can see that, by possibly reducing the value of $t$, the map $x\mapsto (\bar{x},s)$ is a bi-Lipschitz homeomorphism of $(\partial \Omega )^t$ onto $\po\times ]-t,t[$. Accordingly, we shall often use the coordinates $(\bar{x},s)$ to represent the point $x\in (\partial \Omega )^t$.

In this section we consider deformations $\tilde\Omega$ of $\Omega$ of the form
\begin{equation}
\tilde\Omega=(\Omega\setminus (\partial \Omega )^t )\cup \{ (\bar{x},s)\in (\partial \Omega )^t \, : \; s <g(\bar{x})\}
\label{tildeom}
\end{equation}
for appropriate functions $g$ on $\po$.

\begin{definition}
\label{class_normal}
Let $\Omega$ and $t$ be as above. Let $0<\rho<t$ and $M>0$. We say that the domain $\tilde\Omega$ belongs to the class $\cC^{m,1}_M(\Omega,t,\rho)$,
$m=0 \mbox{ or }1$, if $\tilde\Omega$ is given by (\ref{tildeom}) for some $C^{m,1}$ function $g$ on $\po$ which takes values in $]-t+\rho,t[$ and satisfies $|g|_{m,1}\leq M$.
\end{definition}

Given $\Gamma \subset \partial \Omega $ and $\tilde \Omega \in \cC^{m,1}_M(\Omega,t,\rho)$, the set $\tilde \Gamma \subset \partial \tilde \Omega $ where homogeneous Dirichlet boundary conditions are imposed, will be given by
\begin{equation}
\label{tildegamma}
\tilde \Gamma =\{(\bar x, g(\bar x)):\ \bar x\in \Gamma   \}.
\end{equation}

\begin{theorem}
Let $\Omega $ be an open set of class $C^{2}$ and $t >0$ be as above. Let $\Gamma $ be an open subset of
$\partial \Omega $ with Lipschitz continuous boundary in $\partial\Omega  $.
Then there exists $2<q_0\leq\infty$ such that
for any $r>\max \{ 2,N(q_0-1)/q_0\} $
 the following statements hold:
\begin{itemize}
\item[(i)]
There exists $c_1>0$ such that
\begin{equation}
\left(\sum_{n=1}^{\infty}\left| \frac{1}{
{\lambda}_n[\tilde\Omega , \tilde \Gamma ] +1} - \frac{1}{{\lambda}_n[\Omega , \Gamma ] +1}\right|^r\right)^{1/r}
\leq c_1 |\Omega \vartriangle \tilde \Omega  |^{\frac{q_0-2}{rq_0}} ,
\end{equation}
for all $\tilde \Omega \in  \cC^{0,1}_M(\Omega,t,\rho)$  such that,
$|\Omega \vartriangle \tilde \Omega |\le c_1^{-1}
$,
where  $\tilde \Gamma \subset \tilde\Omega$ is given by  (\ref{tildegamma}).
\item[(ii)] Let $\lambda [\Omega ,\Gamma ]$ be an eigenvalue of multiplicity $m$ and let $n\in \N$ be such that
$\lambda [\Omega ,\Gamma ] =\lambda_n[\Omega ,\Gamma ]=\dots =\lambda_{n+m-1}[\Omega , \Gamma ]$.
There exists $c_2>0$ such that the following is true: if
$\tilde \Omega \in   \cC^{0,1}_M(\Omega,t,\rho) $,
$
|\Omega \vartriangle \tilde \Omega |\le c_2^{-1},
$
and  $\tilde \Gamma\subset \partial \tilde \Omega $ is given by (\ref{tildegamma}) then, given orthonormal eigenfunctions
$\psi_n[\tilde \Omega ,\tilde \Gamma ],  \dots ,$ $ \psi_{n+m-1}[\tilde \Omega ,\tilde \Gamma ]$, there exist orthonormal eigenfunctions $ \psi_n[\Omega ,\Gamma ], \dots ,$ $ \psi_{n+m-1}[\Omega ,\Gamma ] $ such that
$$
\| \psi_n[\Omega ,\Gamma ]-\psi_n[\tilde \Omega,\tilde \Gamma ] \|_{L^2(\Omega \cup\tilde \Omega)}\le c_2|\Omega \vartriangle\tilde \Omega |^{\frac{q_0-2}{rq_0}}.
$$
\end{itemize}
Moreover, if in addition $A_{ij}\in C^{0,1}( \R^N )$, $ \tilde \Omega \in  \cC^{1,1}_M(\Omega,t,\rho) $ and $\Gamma $ is a connected component of
$\partial \Omega $ then statements (i) and (ii) hold with $q_0=\infty $.

\label{thm:C11bis}
\end{theorem}
{\em Proof.} The proof is essentially a repetition of the proof of Theorem \ref{thm:C11}: the transformation $\Phi$
is defined as in Lemma  \ref{graf}, with $\partial \Omega $ replacing $W$ and curvilinear coordinates $(\bar x, s)$ replacing the
local euclidean coordinates $(\bar x , x_N)$. \hfill $\Box$

\section{Appendix}

In this section we briefly discuss how  Theorem~\ref{mainthm} can be used to obtain stability estimates for the solutions of the Poisson problem.

\begin{theorem}
\label{thm:appendix}
Let {\rm (A)} be satisfied. Let the operators $L, \tilde  L, \hat L$ satisfy {\rm (P)}  and  $\Omega $ satisfy the interior cone condition.  Let $f\in L^2(\R^N)$  and let $v\in \cW, \tilde v\in\widetilde \cW$ be such that
\begin{equation}\left\{
\begin{array}{ll}
(L+1)v=f,\ \ & {\rm in}\ \phi (\Omega )  ,\vspace{1mm}\\
(\tilde L +1)\tilde v = f,\ \ & {\rm in}\ \tilde \phi (\Omega ) .
\end{array}\right.
\end{equation}
Let $s=[q_0/(q_0-2)] \max\{ 2, \alpha+2\gamma\}$.
If $N\geq 3$, then
there exists $c>0$ depending only on $N, \tau , \alpha ,c^*, q_0 , C, \gamma ,  \Omega $
such that
\begin{equation}\nonumber
\|v-\tilde v \|_{L^2(\phi (\Omega )\cup \tilde \phi (\Omega ))}\le c\left( ( |\cD |^{1/N}  +\delta_s(\phi ,\tilde\phi)) \| f  \|_{L^2(\R^N)} +\| f\circ \phi-f\circ\tilde\phi \|_{L^2(\Omega )}  \right),
\end{equation}
where
$\cD=\{x\in \Omega :\, \phi (x)\ne \tilde\phi (x)  \}$.
The same is true if $N=2$ provided $|\cD |^{1/N} $ is replaced by $|\cD |^{\frac{1}{2}-\epsilon}$, $\epsilon >0$.
\end{theorem}
{\em Proof.} Observe that
\begin{equation}\left\{
\begin{array}{ll}
(H+1)(v\circ \phi )=f\circ \phi,\ \ & {\rm in}\ \Omega  ,\vspace{1mm}\\
(\tilde H +1)(\tilde v\circ \tilde\phi ) = f\circ\tilde \phi,\ \ & {\rm in}\   \Omega  ,
\end{array}\right.
\end{equation}
hence
\begin{equation}
\| v\circ \phi -\tilde v\circ \tilde\phi \|_{L^2(\Omega )}\le \| f\circ \phi -f\circ\tilde\phi \|_{L^2(\Omega )}+
\| (\tilde H+1)^{-1}-( H+1)^{-1} \| \| f\circ \tilde\phi \|_{L^2(\Omega)}.
\end{equation}
By proceeding as in the proof of Theorem~\ref{holder} one can easily see that
\begin{eqnarray}\lefteqn{
\|  v-\tilde v \|_{L^2(\phi (\Omega )\cup \tilde\phi (\Omega ))} }  \nonumber \\   & & \le  c \left( \| v\circ  \phi -v\circ\tilde \phi    \|_{L^2(\Omega ) }  +
\| \tilde v\circ  \phi -\tilde v\circ\tilde \phi    \|_{L^2(\Omega ) }
 \right. \nonumber \\  & & \left.  + \| f\circ \phi -f\circ\tilde\phi \|_{L^2(\Omega )}+
\| (\tilde H+1)^{-1}-( H+1)^{-1} \| \| f\circ \tilde\phi \|_{L^2(\Omega )}.
 \right)
\end{eqnarray}
By the Sobolev embedding it follows that if $N\geq 3$
\begin{eqnarray*}\lefteqn{
\| v\circ  \phi -v\circ\tilde \phi    \|_{L^2(\Omega ) }, \|\tilde  v\circ  \phi -\tilde v\circ\tilde \phi    \|_{L^2(\Omega ) }}\\ &  &\qquad \quad\le  c
|\cD|^{1/N} (\| u\|_{L^2(\Omega )} +  \| \nabla u \|_{L^2(\Omega )})  \le   c
|\cD|^{1/N} \| f\|_{L^2(\R^N )}.
\end{eqnarray*}
The same is true for $N=2$ provided  $|\cD |^{1/N} $ is replaced by $|\cD |^{\frac{1}{2}-\epsilon}$, $\epsilon >0$.
Moreover by Theorem~\ref{mainthm}  it follows that
$$
\| (\tilde H+1)^{-1}-( H+1)^{-1} \| \| f\circ \tilde\phi \|_{L^2(\Omega )}\le c \delta_s(\phi , \tilde\phi)\| f\|_{L^2(\R^N )} .
$$
Thus, the statement follows by combining the estimates above. \hfill $\Box $

\

We now apply the previous theorem in order to  estimate $ \| u -\tilde u\|_{L^{2}(\Omega\cup \tilde \Omega )}$
where $u$, $\tilde u$ are the solutions to the following mixed boundary valued problems and $\tilde \Omega $ is either a local perturbation of $\Omega $ as in Section~\ref{localsec} or a global normal perturbation as in Section~\ref{globsec}:

\begin{equation}
\label{mixedpoi}
\left\{
\begin{array}{ll}
-\sum_{i,j=1}^N\frac{\partial}{\partial x_i }\big( A_{ij}(x)\frac{\partial u}{\partial x_j} \big)=f, &\ {\rm in}\ \Omega ,
\vspace{0.2cm}\\
u=0, &\ {\rm on}\ \Gamma ,\vspace{0.2cm}\\
\sum_{i,j=1}^N A_{ij}\frac{\partial u}{\partial x_j}\nu_{i}=0, &\ {\rm on}\ \partial\Omega  \setminus \Gamma ,
\end{array}\right.
\end{equation}
\vspace{12pt}
\begin{equation}
\label{mixedpoitilde}
\left\{
\begin{array}{ll}
-\sum_{i,j=1}^N\frac{\partial}{\partial x_i }\big( A_{ij}(x)\frac{\partial \tilde u}{\partial x_j} \big)=f, &\ {\rm in}\ \tilde \Omega ,
\vspace{0.2cm}\\
\tilde u=0, &\ {\rm on}\ \tilde \Gamma ,\vspace{0.2cm}\\
\sum_{i,j=1}^N A_{ij}\frac{\partial \tilde u}{\partial x_j}\nu_{i}=0, &\ {\rm on}\ \partial\tilde \Omega  \setminus \tilde \Gamma .
\end{array}\right.
\end{equation}

For any $s>0$ we set
$$
{\mathcal{M}}_f(s )=\sup_{\substack{A\subset \R^N\\ |A|\le s}}\left(\int_{A}|f|^2dx\right)^{1/2}.
$$
The next theorem is a simple corollary of Theorem \ref{thm:appendix} and inequality \eqref{meas}.
\begin{theorem} Let $\Omega , \tilde \Omega , \Gamma ,\tilde \Gamma$ be either as in Theorem \ref{thm:C11} or as
in Theorem \ref{thm:C11bis}. Then the following is true: there exists $2<q_0\le \infty $ such that for any  $r>\max \{ 2,N(q_0-1)/q_0\} $ there exists $c>0$ such that
if $|\Omega\triangle \tilde\Omega|<c^{-1}$ then
\begin{equation}
\label{lastest}
\| u-\tilde u \|_{L^2(\Omega \cup\tilde \Omega)}\le c\left(|\Omega \vartriangle\tilde \Omega |^{\frac{q_0-2}{rq_0}}\| f\|_{L^2(\R^N)} + {\mathcal{M}}_f(c |\Omega \vartriangle\tilde \Omega | )\right).
\end{equation}
Moreover, if in addition $A_{ij}\in C^{0,1}(\R^N)$, $\Omega, \tilde \Omega \in  C^{1,1}$ and $\Gamma$ is a connected component of
$\partial \Omega $ then estimate (\ref{lastest}) holds with $q_0=\infty $.
\end{theorem}

{\bf Acknowledgments.} This work was supported by the research project `Problemi di stabilit\`{a} per operatori differenziali' of the University of Padova, Italy. The third author expresses his gratitude to the  Department of Mathematics of the University of Athens for the kind hospitality during the preparation of this paper.

\vspace{1cm}\noindent
{\it Gerassimos Barbatis\\ Department of Mathematics\\ University of Athens\\ 157 84 Athens\\ Greece\\ gbarbatis@math.uoa.gr\vspace{2mm}\\
Victor~I. Burenkov\\ Pier Domenico Lamberti\\
Dipartimento di Matematica Pura ed Applicata\\ Universit\`{a} degli Studi di Padova\\ Via Trieste, 63\\ 35121 Padova\\ Italy\\ burenkov@math.unipd.it\\ lamberti@math.unipd.it


\begin{thebibliography}{RRR}

\bibitem{bar1}
G.~Barbatis, Spectral stability under $L\sp p$-perturbation of the second-order coefficients, {\it   J. Differential Equations}, {\bf 124},  1996,   pp.~302-323.

\bibitem{bar3}
G.~Barbatis, Stability and regularity of higher order elliptic operators with measurable coefficients, {\it J. London Math. Soc.}, {\bf  58}, 1998, pp.~342-352.

\bibitem{buda} V.I.~Burenkov, E.B.~Davies, Spectral stability of the Neumann
Laplacian,  {\it J.~Differential Equations}, {\bf 186}, pp.~485-508, 2002

\bibitem{bulahigh} V.I.~Burenkov, P.D.~Lamberti, Spectral stability of higher order uniformly elliptic operators, in Sobolev Spaces in Mathematics II.
Applications in Analysis and Partial Differential Equations (to the centenary of Sergey Sobolev), edited by V. 	Maz'ya, {\it International Mathematical Series}, {\bf 9},  Springer, New York, 2008, pp.~69-102.

\bibitem{bula} V.I.~Burenkov, P.D.~Lamberti, {Spectral stability  of Dirichlet second order uniformly elliptic operators}, {\it J.~Differential Equations}, {\bf 244},
pp.~1712-1740, 2008.

\bibitem{bulaneu} V.I.~Burenkov, P.D.~Lamberti, {Spectral stability of general non-negative self-adjoint operators with applications to Neumann-type operators.}, {\it J.~Differential Equations}, {\bf 233},
pp.~345-379, 2007.

\bibitem{bulalanz}  V.I.~Burenkov, P.D.~Lamberti, M.~Lanza de Cristoforis, Spectral stability of nonnegative selfadjoint operators. (Russian)  {\it  Sovrem. Mat. Fundam. Napravl.},   {\bf 15},   2006, pp.~76-111. English translation in: {\it Journal of Mathematical Sciences}, {\bf  149},  2008, pp.~1417-1452.

\bibitem{bulanz}  V.I.~Burenkov, M. Lanza de Cristoforis,  Spectral Stability of the Robin Laplacian,  {\em Proceedings of the Steklov
Institute of Mathematics}, {\bf 260}, 2008, pp.~68-89.

\bibitem{da2000} E.B.~Davies, Sharp boundary estimates for elliptic operators., {\it  Math. Proc. Cambridge Philos. Soc.},  {\bf  129},  2000, pp.~165--178.

\bibitem{da1993}  E.B.~Davies, Eigenvalue stability bounds via weighted Sobolev spaces., {\it   Math. Z.},  {\bf 214}, 1993, pp.~357-371.

\bibitem{daheat} E.B.~Davies, Heat kernels and spectral theory , Cambridge University Press, Cambridge, 1989.

\bibitem{D}P.~Deift, {Applications of a commutation formula}, {\it Duke Math. J.}, {\bf 45}, 1978, pp.~267-309

\bibitem{groger} K.~Gr\"{o}ger,
{A $W\sp {1,p}$-estimate for solutions to mixed boundary value problems for second order elliptic differential equations},
{\em Math. Ann.}, {\bf 283}, 1989, pp.~679-687.


\bibitem{he}
D.~Henry, {Perturbation of the boundary in boundary-value problems of partial differential equations}, London Mathematical Society Lecture Note Series, 318, Cambridge University Press, Cambridge, 2005.
\bibitem{m5}
V.A. Kozlov, V.G. Maz'ya, A.B. Movchan, Asymptotic analysis of fields in multi-structures, Oxford Mathematical Monographs. Oxford Science Publications. The Clarendon Press, Oxford University Press, New York, 1999. xvi+282 pp.
\bibitem{lala03}
P.D.~Lamberti, M.~Lanza de Cristoforis, A global Lipschitz
continuity result for a domain dependent Dirichlet eigenvalue
problem for the Laplace operator, {\it Z. Anal. Anwendungen},
{\bf 24}, 2005, pp.~277-304.

\bibitem{lala03neu}
P.D.~Lamberti, M.~Lanza de Cristoforis,
A global Lipschitz continuity result for a domain-dependent Neumann eigenvalue problem for the Laplace operator, {\it  J. Differential Equations}, {\bf  216},  2005,  pp.~109-133.


\bibitem{m1} V. Maz'ya, S. Nazarov, B. Plamenevskii, Asymptotic theory of elliptic boundary value problems in singularly perturbed domains. Vol. I. 
{\it Operator Theory: Advances and Applications,} 111. Birkh\"{a}user Verlag, Basel, 2000. xxiv+435 pp. 


\bibitem{m2} V. Maz'ya, S. Nazarov, B. Plamenevskii, Asymptotic theory of elliptic boundary value problems in singularly perturbed domains. Vol. II.
{\it Operator Theory: Advances and Applications,} 112. Birkh\"{a}user Verlag, Basel, 2000. xxiv+323 pp. 





{\it Asymptot. Anal.}, {\bf 52}, 2007, pp.~173--206. 



\bibitem{mey} N.G.~Meyers, An $L\sp{p}$-estimate for the gradient of solutions of second order elliptic divergence equations, {\it Ann. Scuola Norm Sup.
Pisa}, {\bf 17},  1963, pp.~189-206.

\bibitem{netsa} Y.~Netrusov, Y.~Safarov, Weyl asymptotic formula for the Laplacian on domains with rough boundaries, {\it  Comm. Math. Phys.}, {\bf  253},   2005,   pp.~481-509.

\bibitem{Pa1} M.M.H.~Pang,  Approximation of ground state eigenvalues and eigenfunctions of Dirichlet Laplacians. {\it  Bull. London Math. Soc.}, {\bf 29},  1997,  pp.~720-730.


\bibitem{RS} M.~ Reed, B.~Simon, Methods of modern mathematical physics, II Fourier Analysis, Self-adjointness. Academic Press 1975.

\bibitem{sava}  G.~Savar\'{e}, G.~Schimperna,
Domain perturbations and estimates for the solutions of second order elliptic equations.
{\it J. Math. Pures Appl. (9)},  {\bf 81},  2002, pp.~1071-1112.


\bibitem{S}B.~Simon, Trace ideals and their applications, Cambridge University Press, Cambridge-New York, 1979.


\bibitem{tro}
G.M.~Troianiello, {Elliptic Differential Equations and Obstacle
  Problems}, Plenum Press, New York, 1987.

\end{thebibliography}
\end{document}